\documentclass[10pt]{article}

\usepackage{amsfonts}
\usepackage{amssymb}
\usepackage{amsmath}
\usepackage{url}
\usepackage{pdfpages}
\usepackage{indentfirst}
\usepackage{color}

\numberwithin{equation}{section}
\numberwithin{equation}{subsection}

\numberwithin{figure}{section}
\numberwithin{figure}{subsection}

\definecolor{myblue}{rgb}{0,0,1}

\newcommand{\abs}[1]{|#1|}

\newcommand{\w}{{\mathrm w}}
\newcommand{\W}{{\mathrm W}}
\newcommand{\z}{{\mathrm z}}
\newcommand{\RH}{{\mathrm{RH}}}

\newcommand{\myd}{{\mathrm d}}

\newcommand{\Arr}{{\mathrm {Arr}}}
\newcommand{\Bow}{{\mathrm {Bow}}}
\newcommand{\bow}{{\mathrm {bow}}}
\newcommand{\arr}{{\mathrm {arr}}}
\newcommand{\lS}{{\mathrm{Spec}}^\lambda}
\newcommand{\Targ}{{\mathrm{Targ}}}
\newcommand{\targ}{{\mathrm{targ}}}
\newcommand{\Rend}{{\mathrm{Rend}}}
\newcommand{\Orb}{{\mathrm{Orb}}}
\newcommand{\orb}{{\mathrm{orb}}}

\newcommand{\zetastar}{{\zeta^*}}
\newcommand{\zetatilde}{{\tilde\zeta}}
\newcommand{\zetatriv}{\zeta_{\mathrm T}}
\newcommand{\zetahat}{{\hat{\zeta}}}

\author{{\sc Yu.\,V.\,Matiyasevich}}

\title{{
\bf Hidden Life of \\Riemann's Zeta Function}\\
 1. Arrow, Bow, and Targets}

\date{Steklov Institute of Mathematics at St.Petersburg, Russia\\[2mm]
\url{http://logic.pdmi.ras.ru/~yumat}}

\begin{document}

\maketitle

\begin{abstract}
The Riemann Hypothesis is reformulated as statements about the
eigenvalues of certain matrices entries of which are defined
via the Taylor series coefficients of the zeta function. These eigenvalues
demonstrate interesting visual patterns allowing one to state a number
of conjectures.
\end{abstract}

\section*{}
\setcounter{section}{1}
\subsection{The Hypothesis}

One of the most interesting and important objects in number theory is
Riemann's zeta function $\zeta(z)$.
It can be defined for $\Re(z)>1$
by the Dirichlet series
\begin{equation}
\zeta(z)=\sum_{n=1}^{\infty}\frac{1}{n^z}.
\label{zeta}
\end{equation}
The function can be extended to the entire complex $z$-plane
with the exception of the point $z=1$ which is the only pole
of $\zeta(z)$. Points $\z_1=-2,\ \z_2=-4,\ \dots,\ \z_n=-2n,\dots$
are known as the \emph{trivial zeros} of the function $\zeta(z)$.
We have the famous

\

{\bf Riemann Hypothesis (version 1).} \emph{All non-trivial
zeros of the function $\zeta(z)$ lie on the
\emph{critical line} $\Re(z)=\frac{1}{2}$.}

\subsection{Trivial zeroes }

There is a tradition (taking its origin in Riemann's seminal
paper \cite{Riemann}) to get rid of the trivial zeros by
dealing with the function
\begin{equation}\xi(z)=
\pi^{-\frac{z}{2}}
(z-1)\zeta(z)\Gamma(1+\tfrac{z}{2})
\label{xi}\end{equation}
rather than with the function $\zeta(z)$ itself
(we  use  modern notation for this function,
Riemann used $\xi(t)$ to denote the function which is
today denoted $\Xi(t)$).
The poles of the factor $\Gamma(1+\tfrac{z}{2})$
in \eqref{xi} cancel the trivial zeros of $\zeta(z)$ and similarly
the factor $z-1$ cancels the pole of $\zeta(z)$.
The factor $\pi^{-\frac{z}{2}}$ influence neither
zeros nor poles but it allows one to state
the \emph{functional equation} in a pretty form:
\begin{equation}
\xi(z)=\xi(1-z).
\label{funceq}\end{equation}

In this paper we won't deprive zeta function of its
trivial zeros but try to take advantage of our
knowledge of precise positions of these zeros. To this end
we will work with the entire function
\begin{equation}
  \zetastar(z)=2(z-1)\zeta(z).
\label{zetastar}\end{equation}
For our purpose we could also omit the factor
$z-1$ and/or use the factor $\pi^{-\frac{z}{2}}$;
this would change the picture(s) so probably separate paper(s)
will be devoted to these variations. The factor $2$ in \eqref{zetastar}
results in the equality
\begin{equation}
   \zetastar(0)=1
\label{zetastar0}\end{equation}
which slightly simplifies some forthcoming formulas.

According to \eqref{funceq} the non-trivial zeros of $\zeta(z)$  lie
symmetrically around the critical line, so we have

\

{\bf Riemann Hypothesis (version 2).} \emph{The trivial
zeros $\z_1=-2,\ \z_2=-4,\ \dots,\ \z_n=-2n,\dots$
are the only zeros
of the function $\zetastar(z)$ lying in the half-plane
 $\Re(z)<\frac{1}{2}$.}

\subsection{Change of the variable}

A half-plane is a natural object when one deals with Dirichlet
series. However, we are going to deal with Taylor series,
and for them circles are more natural regions. So we make a change
of variable:
\begin{equation}
z=\frac{w}{w+1},\qquad w=\frac{z}{1-z}.
\end{equation}
Under this transformation the critical line becomes
the \emph{critical circle} $|w|=1$, the half-plane $\Re(z)<\frac{1}{2}$
becomes the interior of this circle, and points
\begin{equation}
\w_1=\frac{\z_1}{1-\z_1}=-\frac{2}{3},\
\dots,\ \w_n=\frac{\z_n}{1-\z_n}=-\frac{2n}{2n+1},\dots
\end{equation}
become the \emph{trivial zeros of the function}
\begin{equation}
\zetatilde(w)=\zetastar(\tfrac{w}{w+1}).
\label{zetatilde}\end{equation}

\begin{figure}[t]
\includegraphics[width=.47\textwidth]{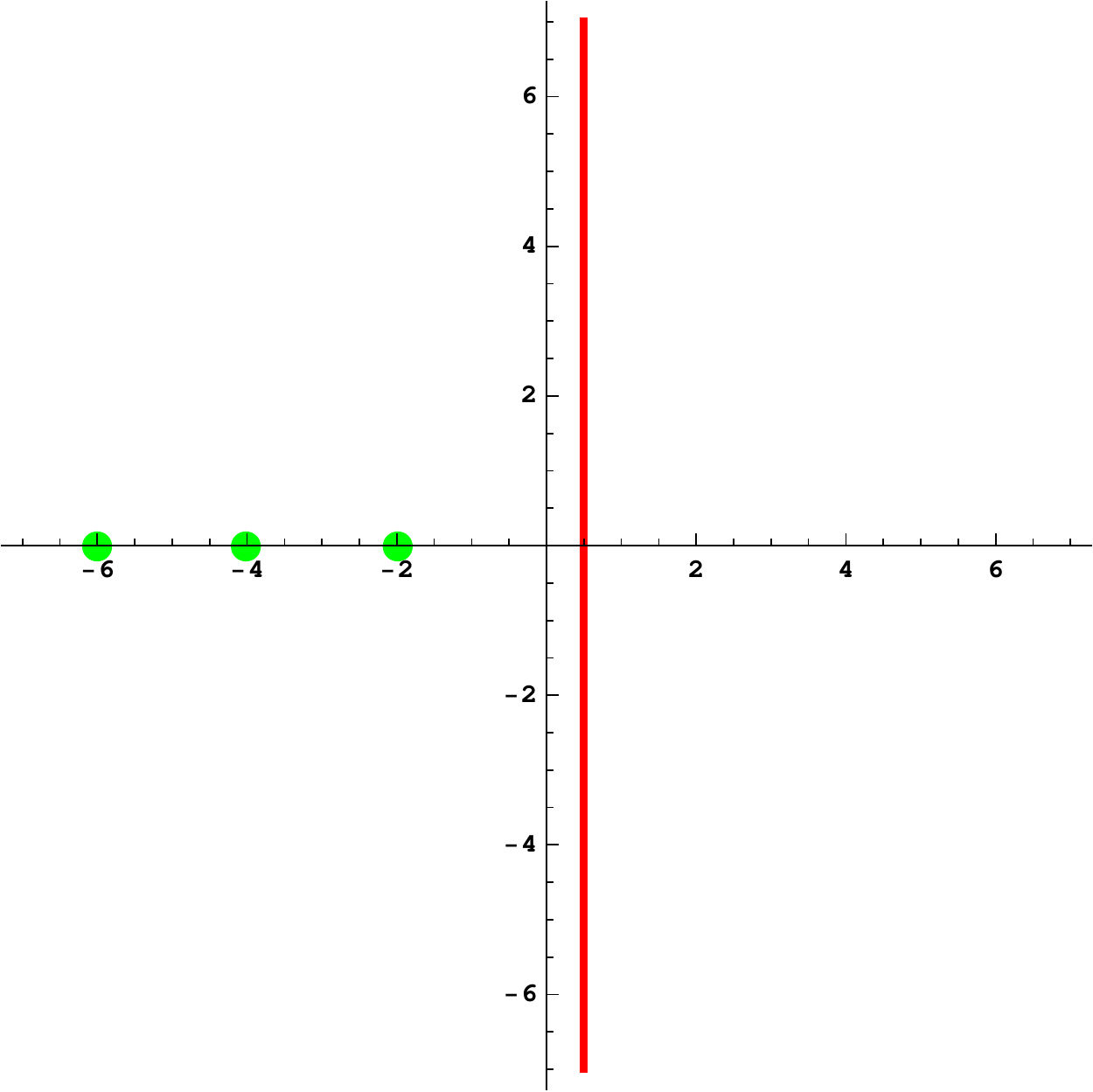}\hfill
\includegraphics[width=.47\textwidth]{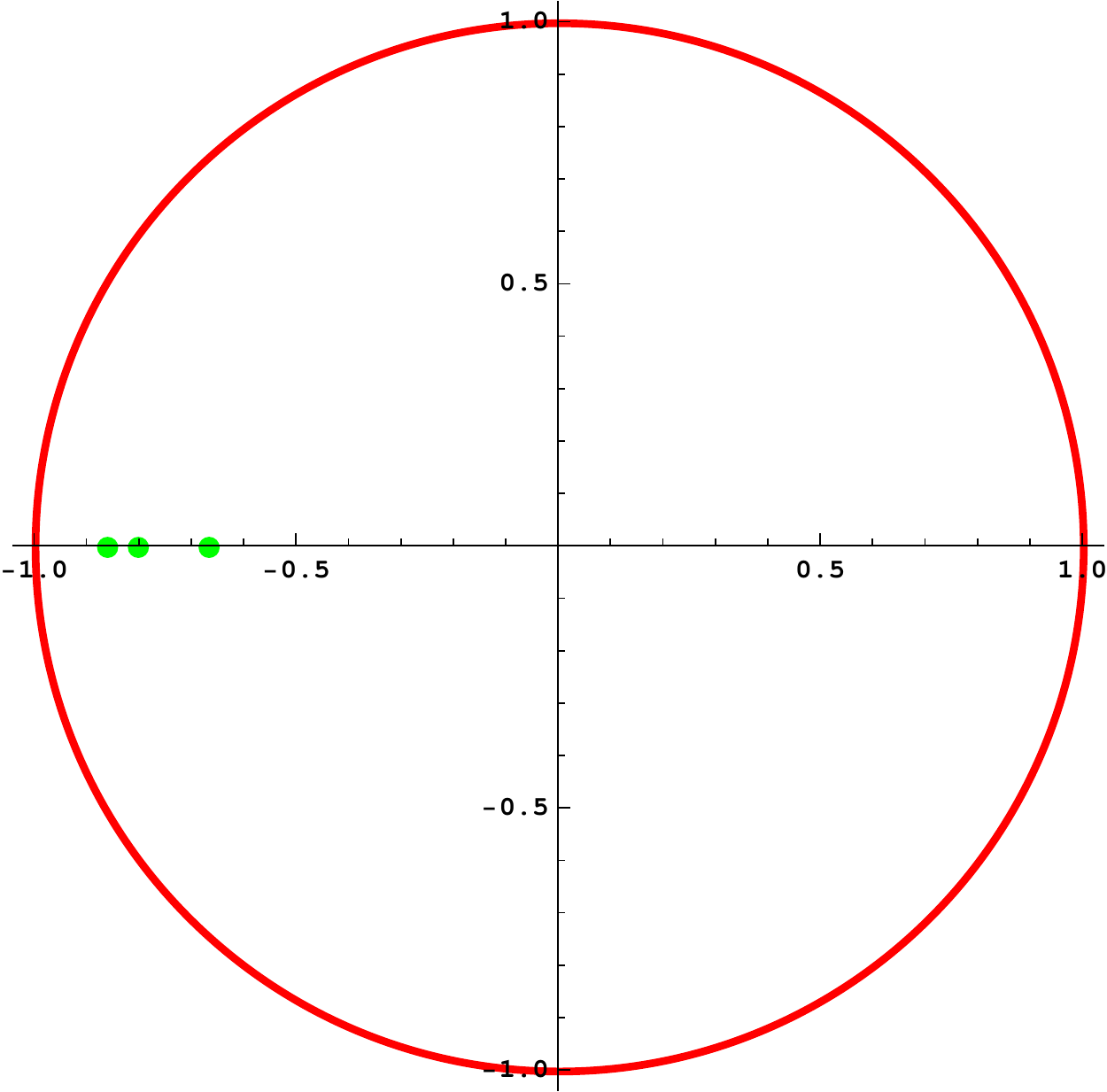}\\
\parbox{.47\textwidth}{\caption{$z$-plane}\label{zplane}}
\hfill\parbox{.47\textwidth}{\caption{$w$-plane}\label{wplane}}
\end{figure}

With this new notation we have

\

{\bf Riemann Hypothesis (version 3).} \emph{The trivial
zeros $\w_1=-\frac{2}{3},\
\dots,\ \w_n=-\frac{2n}{2n+1},\dots$
are the only zeros
of the function $\zetatilde(w)$ lying in the open circle
 $|w|<1$.}

\subsection{Subhypothesises}

It isn't very convenient to work near the critical circle
(full of zeros)
so we split the Riemann  Hypothesis into an infinite series
of weaker statements:

\

{\bf $\mathbf{RH_{\emph{l}}}$, the $\mathbf{\emph{l}}$-th Riemann subhypothesis (version 1).}
 \emph{The trivial
zeros $\w_1=-\frac{2}{3},\
\dots,\ \w_l=-\frac{2l}{2l+1}$
are the only zeros
of the function $\zetatilde(w)$ lying in the closed disk
 $|w|\le \frac{2l+1}{2l+2}$.}

\

While each of the subhypotheses is weaker than the Riemann Hypothesis,
taken together, they are equivalent to it:

\

{\bf Riemann Hypothesis (version 4).} \emph{For every $m$
the subhypothesis $\RH_m$ is true.}

\subsection{First question}

Let us ask a ``na\"\i ve'' question:
\emph{Why is $\RH_1$  true?} More precisely, \emph{how can we see
that $\RH_1$ is true?} To answer this question we can expand the function
$1/\zetatilde(w)$
into the
Taylor series:
\begin{equation}
1/\zetatilde(w)=1+\tau_1w+\dots+\tau_nw^n+\dots
\label{invzetatilde}\end{equation}
Now
 according to $\RH_1$
the point $\w_1=-\frac{2}{3}$ should be the only pole of the
function \eqref{invzetatilde} lying inside the circle $|w|\le \frac{3}{4}$.
Respectively, we have

\

{\bf $\mathbf{RH_1}$ (version 2).}
\emph{For $m\rightarrow \infty$
\begin{equation}
\tau_m\,
=\,\big(-\tfrac{3}{2}\big)^m(R_1+o(1))
\label{RH1ver2}\end{equation}
for some non-zero constant $R_1$.}

\

It is easy to see that
\begin{eqnarray}
R_1&=&\frac{3}{2\zetatilde'(-2/3)}\label{R1tilde}\\
&=&
-\frac{1}{36\zeta'(-2)}\label{R1zeta}\\
&=&0.91228851841347\dots\label{R1num}\\
&>&0
\label{R1pos}\end{eqnarray}
so we have

\

{\bf $\mathbf{RH_1}$ (version 3).}
\emph{
\begin{equation}
\lim_{m\rightarrow \infty}((-1)^m\tau_m)^{\frac{1}{m}}\,
=\,\tfrac{3}{2}.
\label{RH1ver3}\end{equation}
}

\subsection{Determinant representation}

 It is easy to see that coefficients $\tau_1,\ \tau_2,\dots$
from \eqref{invzetatilde} can be expressed in terms of the
 coefficients in the Taylor expansion
 \begin{equation}
\zetatilde(w)=1+\theta_1w+\dots+\theta_mw^m+\dots.
\label{zetatildetaylor}\end{equation}
More precisely,
\begin{equation}
\tau_m=(-1)^m\det(L_{1,m})
\end{equation}
where $L_{1,m}$ is the following Toeplitz matrix:\footnote{Clearly, we can change
the order of the columns in this
matrix and obtain a  Hankel matrix whose determinant has the same
absolute value;
the resulting picture(s) are rather different, they are considered in
\cite{MatiyasevichHidden2};
for our purpose we can also deal with the product of our matrix $L_{1,m}$ and its
Hankel counterpart.}
\begin{equation}
L_{1,m}=
\begin{pmatrix}
\theta_{1}&1&0&\dots&0&0\\
\theta_{2}&\theta_{1}&1&\dots&0&0\\
\theta_{3}&\theta_{2}&\theta_1&\dots&0&0\\
\vdots&\vdots&\vdots&\ddots&\vdots&\vdots\\
\theta_{m-1}&\theta_{m-2}&\theta_{m-3}&\dots&\theta_1&1\\
\theta_{m}&\theta_{m-1}&\theta_{m-2}&\dots&\theta_2&\theta_{1}
\end{pmatrix}.
\label{matrixL1m}
\end{equation}
Then, we have

\

{\bf $\mathbf{RH_1}$ (version 4).}
\emph{For  $m\rightarrow \infty$
\begin{equation}
\det(L _{1,m})\,
=\,\big(-\tfrac{3}{2}\big)^m(R_1+o(1))
\label{det1R1}\end{equation}
with the constant $R_1$ defined by \eqref{R1tilde}--\eqref{R1zeta}}

\

\noindent and

\

{\bf $\mathbf{RH_1}$ (version 5).}
\begin{equation}
\lim_{m\rightarrow\infty}\big(\det(L _{1,m})\big)^\frac{1}{m}
=\,\tfrac{3}{2}.
\end{equation}

\subsection{Eigenvalues on average}

Naturally,
\begin{equation}
\det(L _{1,m})=\lambda_{1,m,1}\lambda_{1,m,2}\dots\lambda_{1,m,m}
\label{lambdas}\end{equation}
where $\lambda_{1,m,1}, \ \lambda_{1,m,2},\ \dots, \ \lambda_{1,m,m}$
are the eigenvalues of the matrix $L_{1,m}$.
Thus, we have

\

{\bf $\mathbf{RH_1}$ (version 6).}
\begin{equation}
\lim_{m\rightarrow \infty}
\left(\prod_{n=1}^m{\lambda_{1,m,n}}\right)^\frac{1}{m}=
\frac{3}{2}.
\label{limlambda32}\end{equation}

The (multi)set $\{\lambda_{1,m,1},\ \lambda_{1,m,2},\dots\lambda_{1,m,m}\}$ will be called
the $\lambda$-\emph{spectrum} and will be denoted~$\lS_{1,m}$.

\subsection{Positions of individual eigenvalues}

According to $\RH_1$ the (geometric) mean of
$\lambda_{1,m,1}$,  $\lambda_{1,m,2}$, \dots,  $\lambda_{1,m,m}$
approaches $\frac{3}{2}$ when $m$ goes to infinity, but neither $\RH_1$
nor $\RH$ itself
 tells us anything directly
about the distribution of these eigenvalues. Are they as random as, say, the imaginary parts
of the non-trivial zeros of $\zeta(z)$?
Do the eigenvalues cluster or are they spread around the whole $w$-plane?
Is there any similarity between
eigenvalues corresponding to different values of $m$?

The author was curious to calculate\footnote{Calculations were done
mainly with {\sc Mathematica} and partly with {\sc Pari} on a personal
computer; larger scale computations are very desirable for getting more
insight.}
 the values of the spectra
$\lS_{1,m}$
for initial values of $m$ and have a look at them. Some pictures are
included in this paper, an updated collection of pictures can
be downloaded from \cite{hidden}.
\begin{figure}[phtb]
\includegraphics[width=.47\textwidth]{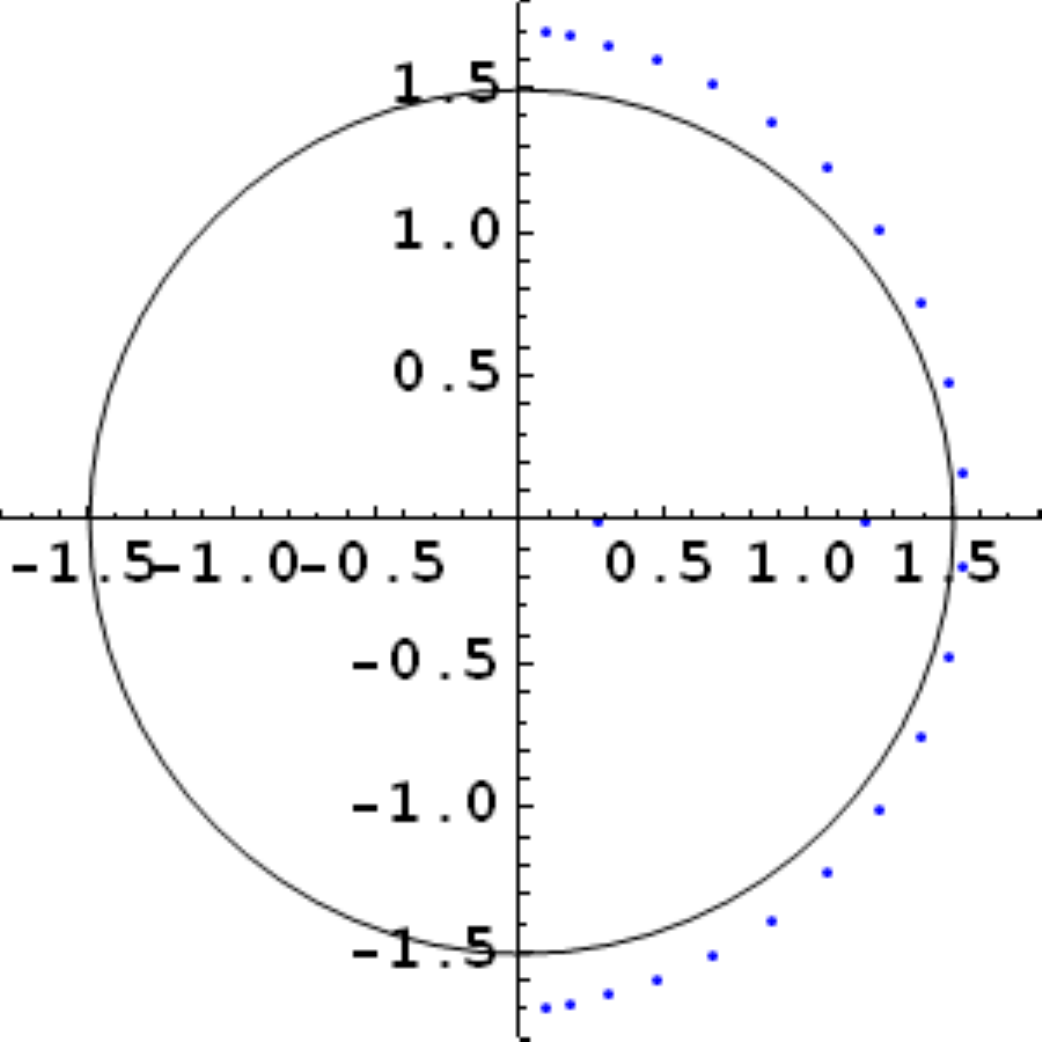}\hfill
\includegraphics[width=.47\textwidth]{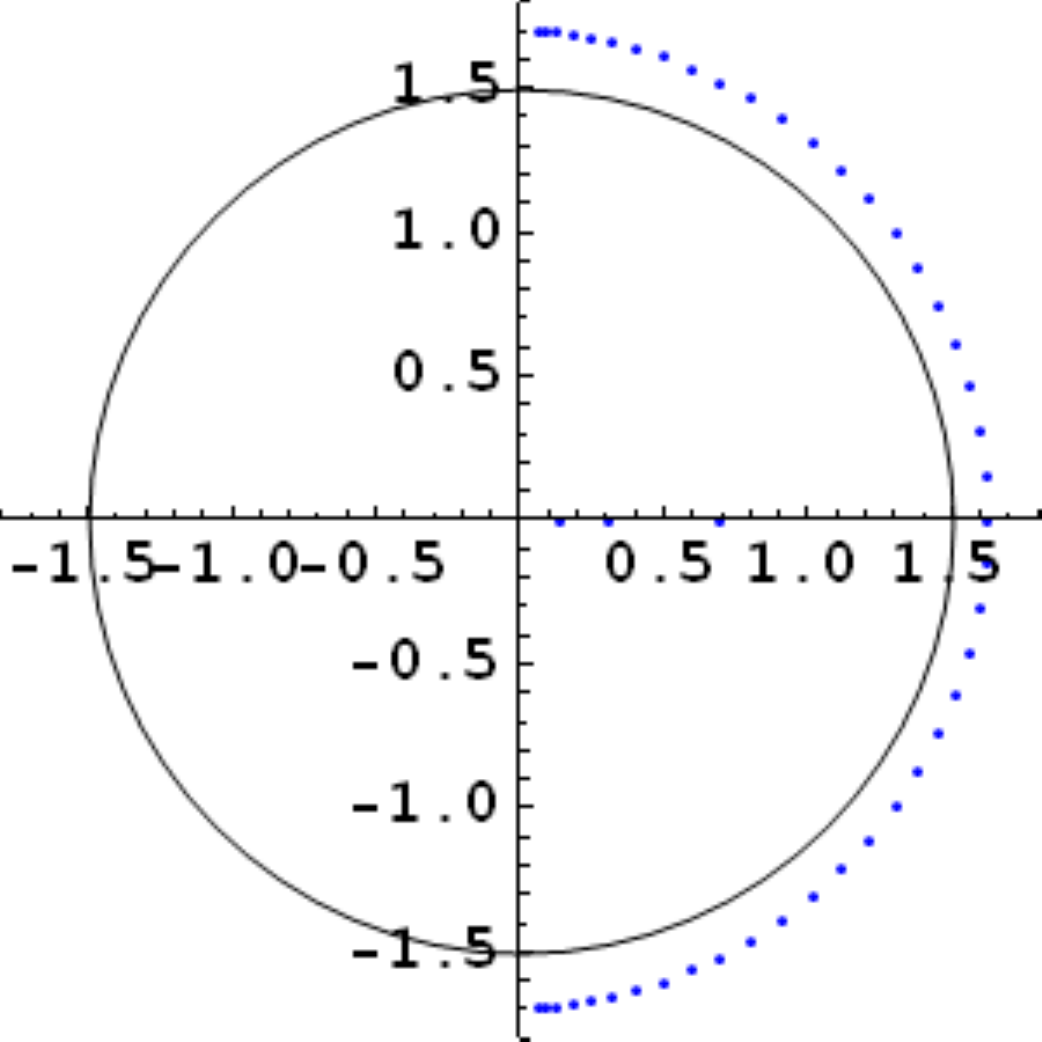}\\
\parbox{.47\textwidth}{\caption{$\lS_{1,24}$}\label{pic124}}
\hfill\parbox{.47\textwidth}{\caption{$\lS_{1,48}$}\label{pic148}}
\end{figure}
\begin{figure}[phtb]
\includegraphics[width=.47\textwidth]{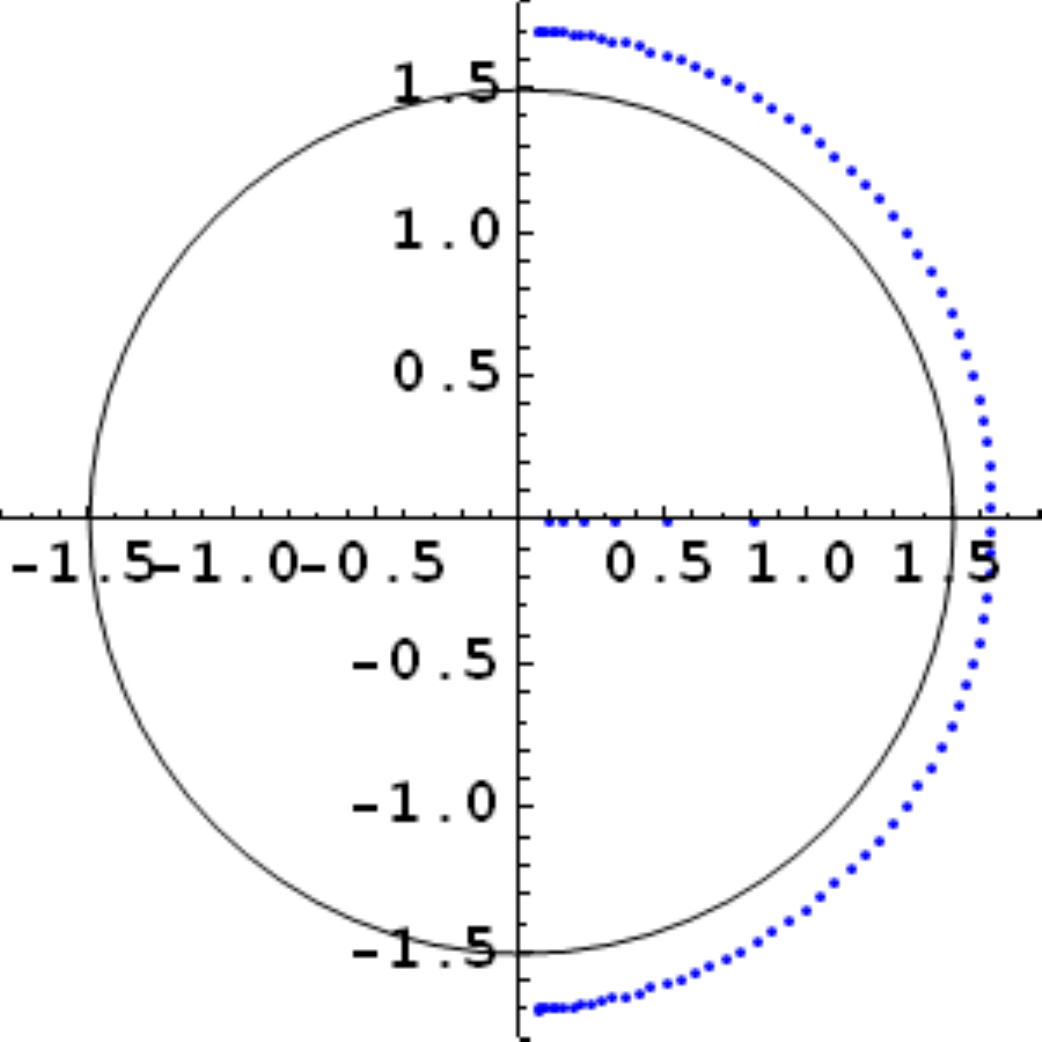}\hfill
\includegraphics[width=.47\textwidth]{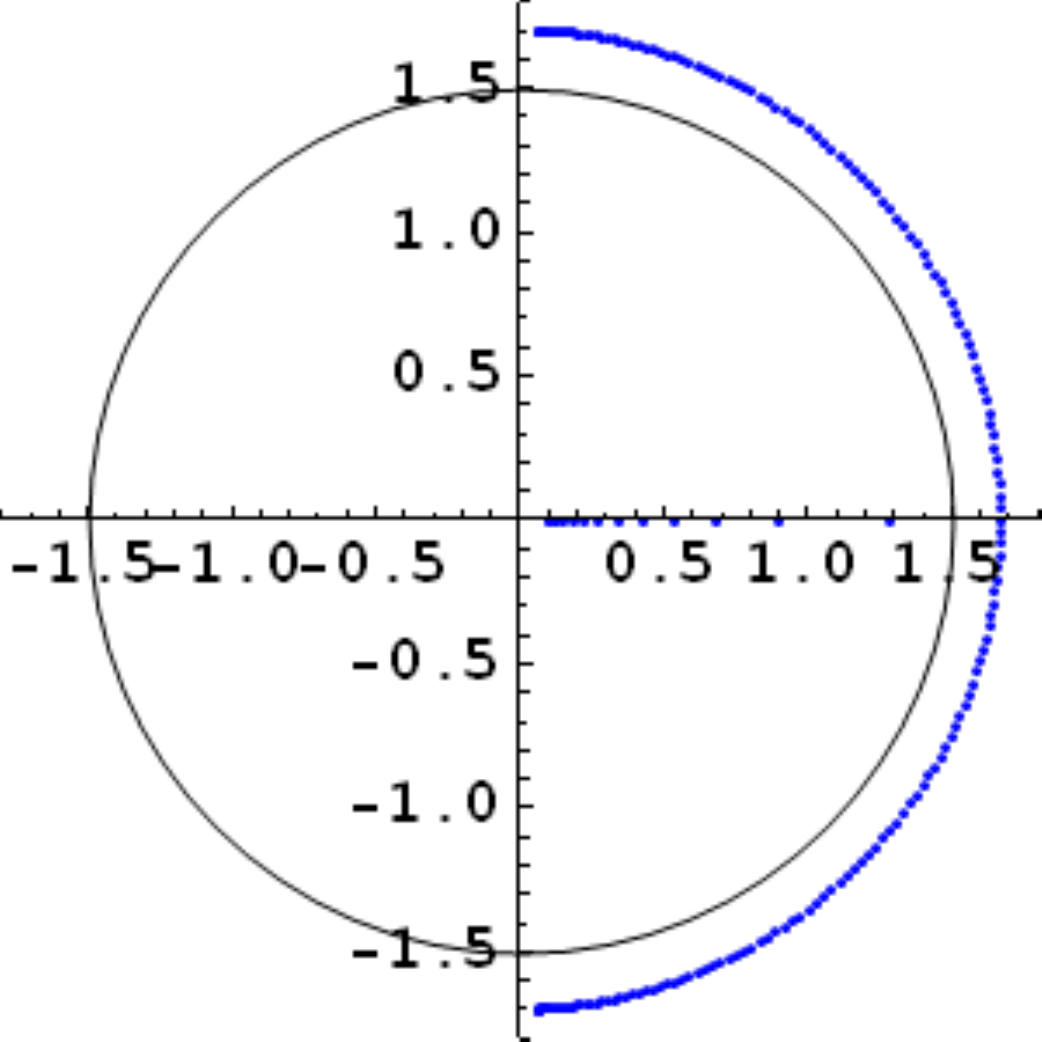}\\
\parbox{.47\textwidth}{\caption{$\lS_{1,96}$}\label{pic196}}\hfill
\parbox{.47\textwidth}{\caption{$\lS_{1,192}$}\label{pic1192}}
\end{figure}
\begin{figure}[tbhp]
\includegraphics[width=\textwidth]{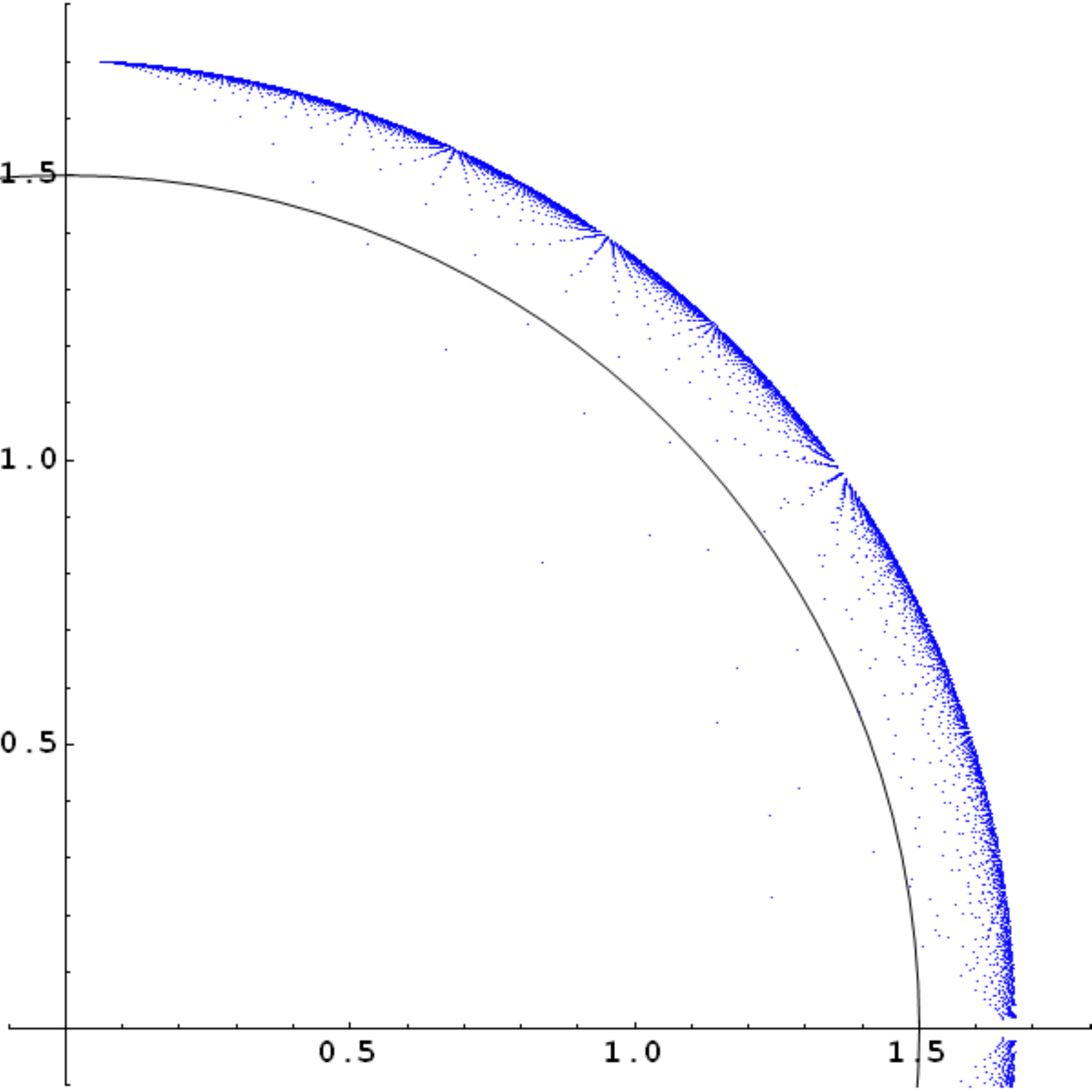}\\
\parbox{\textwidth}{\caption{$\cup_{m=1}^{192}\lS_{1,m}$}\label{pic11192}}
\end{figure}

Figures \ref{pic124}--\ref{pic1192} show  spectra $\lS_{1,m}$
for $m=24,\ 48,\ 96, \ 192$
  respectively together with the circle $\abs{w}=\frac{3}{2}$, the
``ideal'' place for the eigenvalues according to \eqref{limlambda32}.
Figure \ref{pic11192} shows  the union $\cup_{m=1}^{192}\lS_{1,m}$
(in the first quadrant).
The ``hidden life of Riemann's zeta function'' is best seen from an
animation showing the  $\lS_{1,1}$, $\lS_{1,2}$, \dots in succession; such an animation can
be downloaded from \cite{hidden}.

Looking at the pictures we can say that the $\lambda$-spectrum $\lS_{1,m}$
is the union of the  \emph{arrow} $\Arr_{1,m}$
consisting entirely of real eigenvalues
and the   \emph{bow} $\Bow_{1,m}$;
for counting purpose it is reasonable to consider sometimes the
largest real eigenvalue as belonging to the bow rather than to the arrow.
Formally, if the number of real eigenvalues in the spectrum $\lS_{1,m+1}$
is less than the number of real eigenvalues in the spectrum $\lS_{1,m}$, then
we consider the largest real eigenvalue in the spectrum $\lS_{1,m}$
as belonging to $\Bow_{1,m}$ and not belonging to arrow $\Arr_{1,m}$.

\subsection{First Conjectures}

The above pictures suggest the
following conjectures.

\

{\bf Conjecture $\mathbf{1A_1}$.} \emph{There are no multiple eigenvalues.}

\

{\bf Conjecture $\mathbf{1B_1}$.} \emph{$\mathop{\sup}_m( \max (\Arr_{1,m}))$ is a positive number.}

\

{\bf Conjecture $\mathbf{1C_1}$.} \emph{$\inf_m(\min (\Arr_{1,m}))$ is a positive number.}

\

{\bf Conjecture $\mathbf{1D_1}$.} \emph{The numbers $\arr_{1,m}=||\Arr_{1,m}||$
and $\bow_{1,m}=||\Bow_{1,m}||$ of
eigenvalues belonging to
the arrow $\Arr_{1,m}$ and to the bow $\Bow_{1,m}$ respectively
don't decrease when $m$ increases.}

\

{\bf Conjecture $\mathbf{1E_1}$.} \emph{If
\begin{eqnarray}
\Arr_{1,m}&=&\{\lambda_{1,m,1},\lambda_{1,m,2},\dots,\lambda_{1,m,\arr_{1,m}}\},\\
\Arr_{1,m+1}&=&\{\lambda_{1,m+1,1},\lambda_{1,m+1,2},\dots,\lambda_{1,m+1,\arr_{1,m+1}}\}
\end{eqnarray}
and
\begin{eqnarray}
&\lambda_{1,m,1}<\lambda_{1,m,2}<\dots<\lambda_{1,m,\arr_{1,m}},\\
&\lambda_{1,m+1,1}<\lambda_{1,m+1,2}<\dots<\lambda_{1,m+1,\arr_{1,m+1}}
\end{eqnarray}
then
\begin{equation}
\lambda_{1,m+1,1}<\lambda_{1,m,1},\quad
\lambda_{1,m+1,2}<\lambda_{1,m,2},\dots
\lambda_{1,m+1,\arr_{1,m}}<\lambda_{1,m,\arr_{1,m}}.
\end{equation}
}

\

{\bf Conjecture $\mathbf{1F_1}$.}
\emph{Assign the weight $\frac{1}{m}$ to each of the points
$\lambda_{1,m,1},\lambda_{1,m,2},\dots,\lambda_{1,m,m}$
and denote by $\lambda_{1,m}$ corresponding discrete
measure.  Then
\begin{itemize}
\item[$\mathbf{1F'_1.}$] there exists a limiting continuous measure $\lambda_{1}(w)$
concentrated on a ``limiting bow'' and a ``limiting arrow'';
 \item[$\mathbf{1F''_1.}$]
$\int \log(w)\myd\lambda_{1}(w)=\log\left(\frac{3}{2}\right)$.
\end{itemize}
}

\

Clearly, Conjecture $1\mathrm{F}_1$ implies Subhypothesis $\RH_1$.

\subsection{Purely Trivial Zeros}

It is natural to try to understand to what extent the
distribution of
$\lambda_{1,m,1}$, $\lambda_{1,m,2}$, \dots, $\lambda_{1,m,m}$
is due to the trivial zeros, and what is the contribution
of the non-trivial zeros. To this end we can consider the function
\begin{eqnarray}
\zetatriv(z)&=&\frac{\zetastar(z)}{2\xi(z)}\label{zetatriv1} \\
&=&\frac{\pi^{\frac{z}{2}}}
      {\Gamma(1+\frac{z}{2})}.
\label{zetatriv2}\end{eqnarray}
The points
 $\z_1=-2,\ \z_2=-4,\ \dots,\ \z_n=-2n,\dots$
are the only zeros of the function $\zetatriv(z)$.
The factor $2$ in the denominator of \eqref{zetatriv1}
implies the equality
\begin{equation}
\zetatriv(0)=1
\label{zetatriv0}\end{equation}
analogous to the equality \eqref{zetastar0}.

By analogy with \eqref{zetatilde} and \eqref{zetatildetaylor}, for every
analytic function
$f(z)$ such that
\begin{equation}
f(0)=1
\label{f01}\end{equation} we can consider
the
transformed function
\begin{equation}
\tilde f(w)=f\left(\tfrac{w}{1+w}\right)
\end{equation}
with the expansion
 \begin{equation}
\tilde f(w)=1+\theta_1(f)w+\dots+\theta_m(f)w^m+\dots,
\label{ftaylor}\end{equation}
form the matrices $L_{1,m}(f)$, counterparts of \eqref{matrixL1m},
with eigenvalues
$\lambda_{1,m,1}(f)$,  $\lambda_{1,m,2}(f)$ \dots,  $\lambda_{1,m,m}(f)$,
and state various versions of subhypothesis $\RH_1(f)$.
\begin{figure}[p]\centering
\includegraphics[width=.47\textwidth]{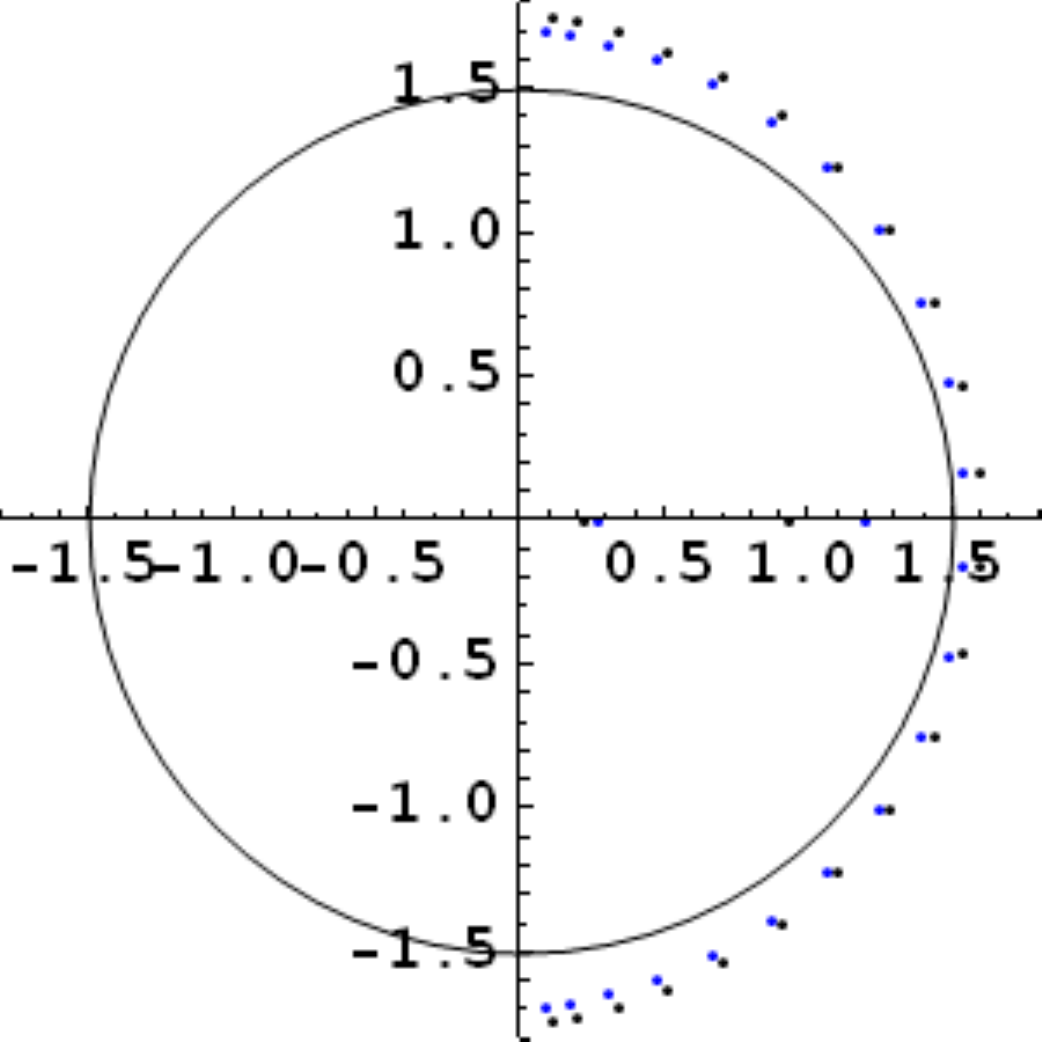}\hfill
\includegraphics[width=.47\textwidth]{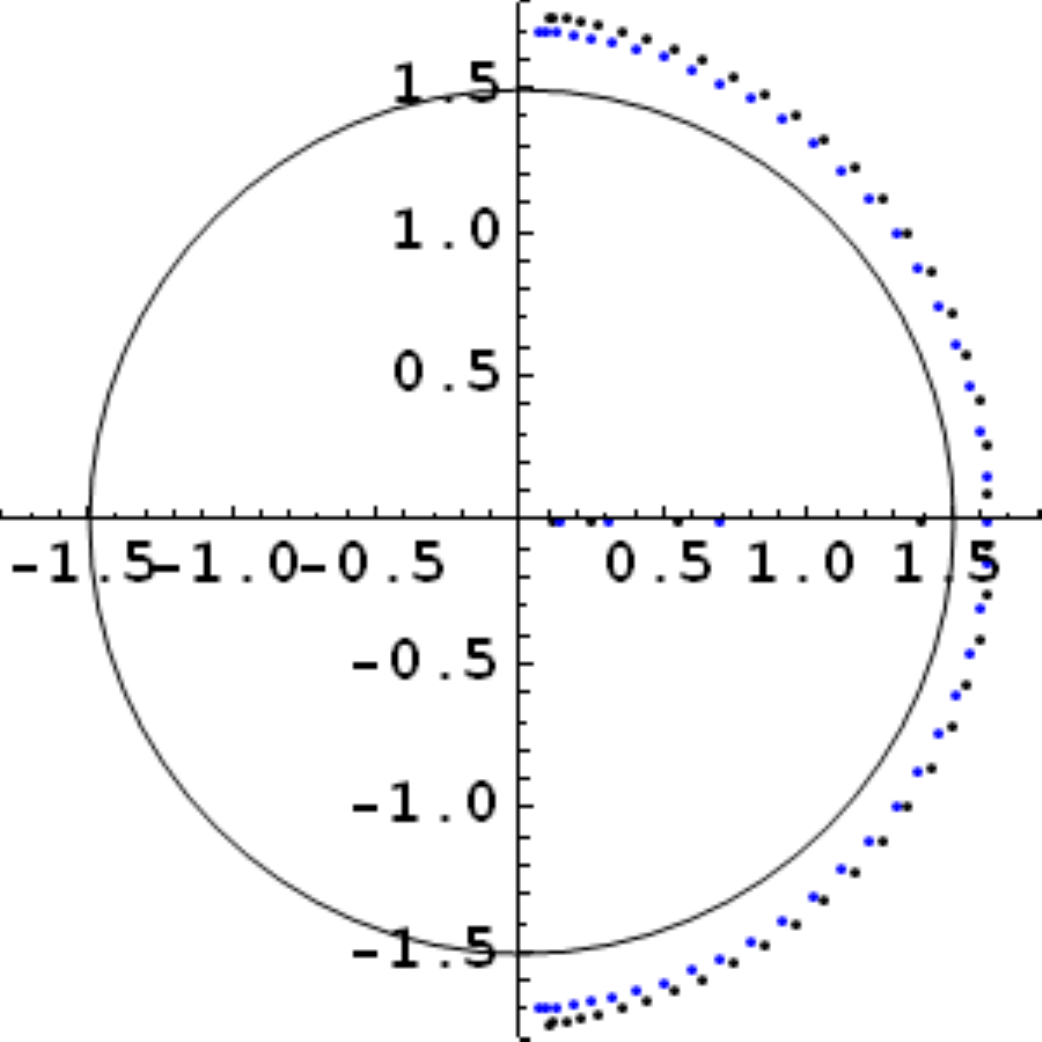}\\
\parbox{.47\textwidth}
{\caption{{\color{myblue}{$\lS_{1,24}(\zetastar)$}} and $\lS_{1,24}(\zetatriv)$\label{pic124T}}}
\hfill\parbox{.47\textwidth}
{\caption{{\color{myblue}{$\lS_{1,48}(\zetastar)$}} and $\lS_{1,48}(\zetatriv)$\label{pic148T}}}
\end{figure}
\begin{figure}[p]
\includegraphics[width=.47\textwidth]{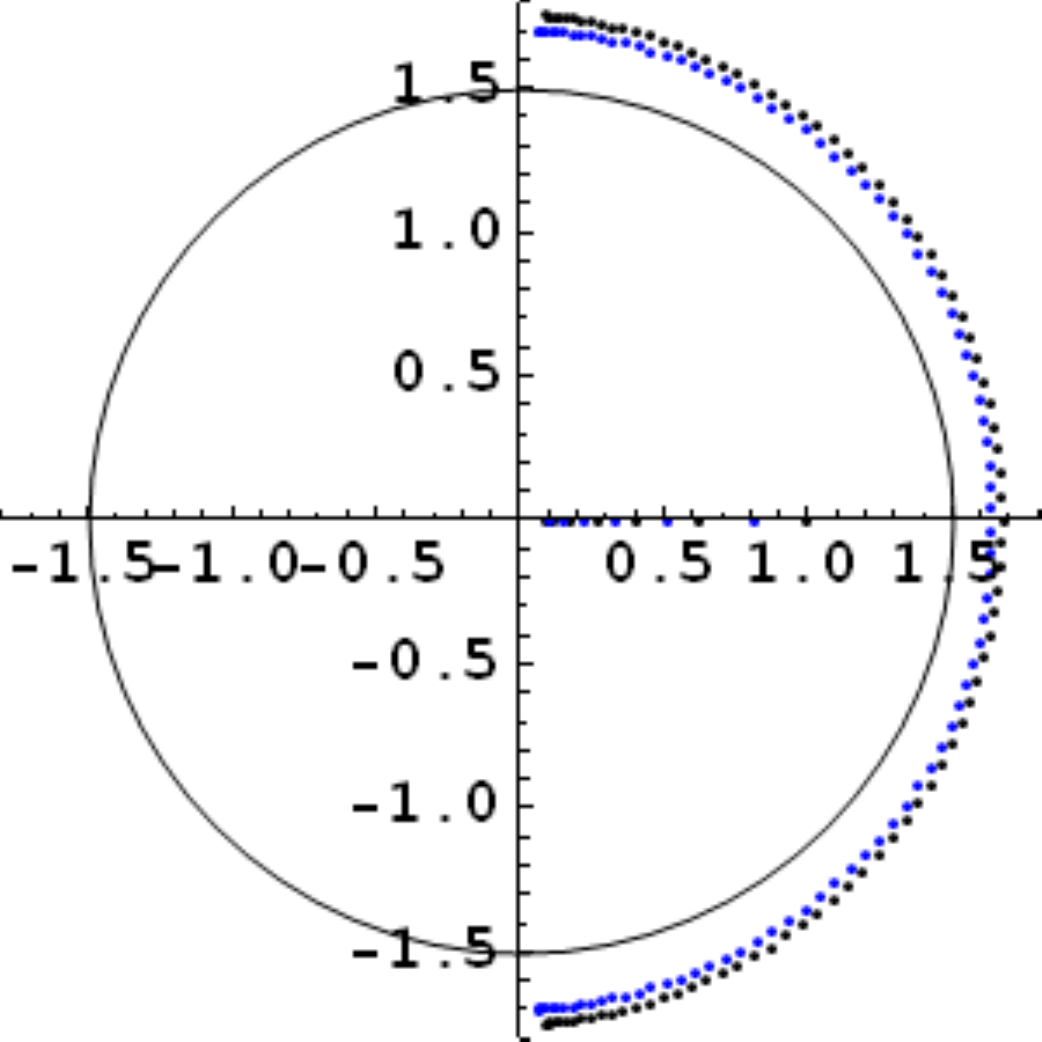}\hfill
\includegraphics[width=.47\textwidth]{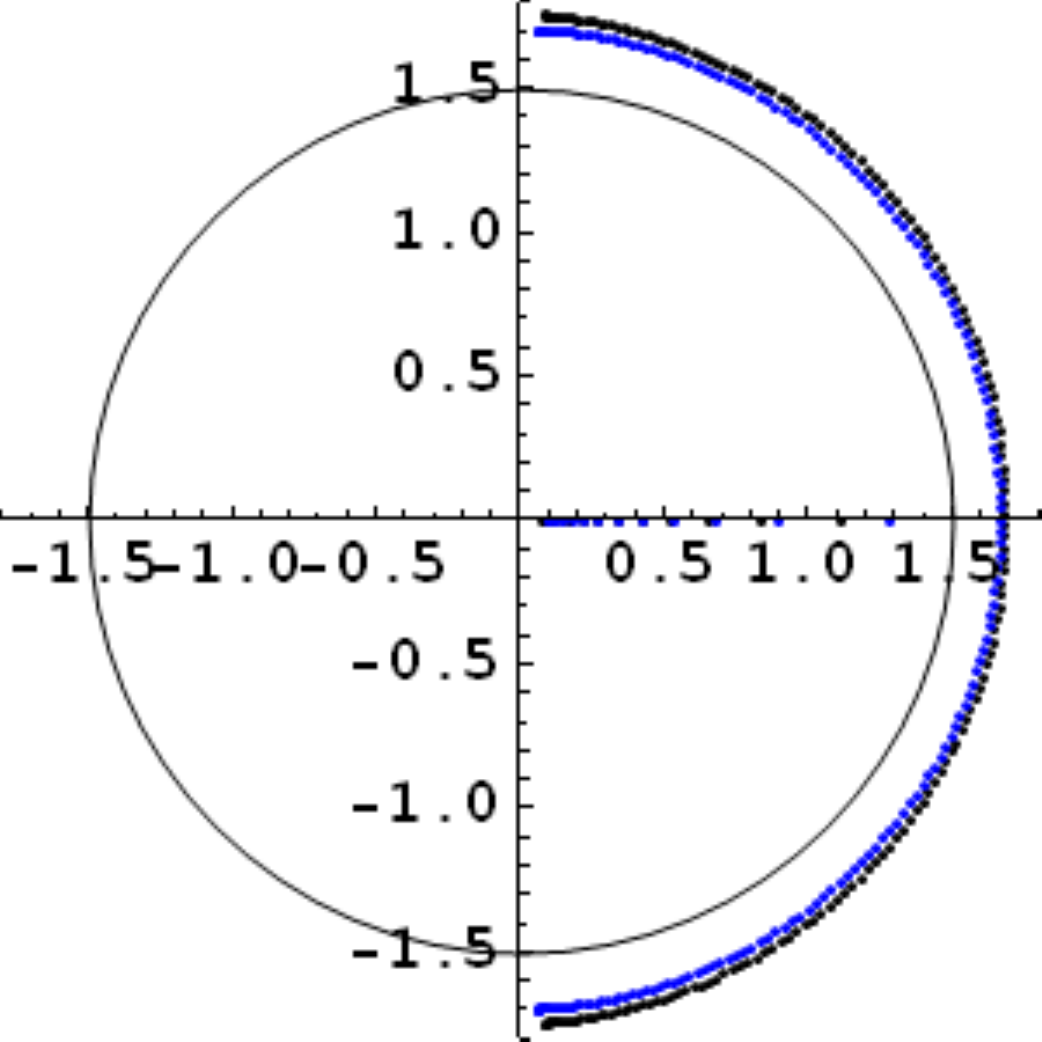}\\
\parbox{.47\textwidth}
{\caption{{\color{myblue}{$\lS_{1,96}(\zetastar)$}} and $\lS_{1,96}(\zetatriv)$\label{pic196T}}}\hfill
\parbox{.47\textwidth}
{\caption{{\color{myblue}{$\lS_{1,192}(\zetastar)$}}  and $\lS_{1,192}(\zetatriv)$\label{pic1192T}}}
\end{figure}

Figures \ref{pic124T}--\ref{pic1192T} show  spectra $\lS_{1,m}(\zetatriv)$ in black color together
with  $\lS_{1,m}(\zetastar)$ in blue  color
for $m=24,\ 48,\ 96, \ 192$
  respectively.
%Figure \ref{pic11192T} shows  the union $\cup_{m=1}^{192}\lS_{1,m}(\zetatriv)$.
These figures suggest that the distribution of the $\lambda$'s is to a great
extent determined by the trivial zeros.

The gamma function is supposed to be ``simple'', ``completely understood'',
a function about which we know everything; it would be natural, as a first step
towards the Riemann Hypothesis, to understand the character of the numbers
$\lambda_{1,m,n}(\zetatriv)$.

\subsection{Further Questions}

Now, how could we see that $\RH_2$, $\RH_3$, \dots are true?
The Taylor expansion
\eqref{invzetatilde} doesn't tell us anything directly about the other
 poles of the function $1/\zetatilde(w)$. One
way to overcome this obstacle could be to consider the function
\begin{equation}
 \zetahat_l(z)=\frac{\zetastar(z)}{\prod_{k=1}^{l-1}(1-\frac{z}{\z_k})};
\end{equation}
a separate paper may be devoted to the corresponding eigenvalues
$\lambda_{1,m,1}(\zetahat_l)$,  $\lambda_{1,m,2}(\zetahat_l)$ \dots,
$\lambda_{1,m,m}(\zetahat_l)$.

\subsection{Pad{\'e} approximations}

However, another approach to $\RH_m$ for arbitrary $m$ looks more promising.
This approach treats the first $m$ trivial zeros ``on equal'' and is
based on Pad{\'e} approximations.

To begin with, let $P_{1,m}(w)$ and $Q_{1,m}(w)$ be polynomials such that
\begin{eqnarray}
\zetatilde(w)\approx\frac{P_{1,m}(w)}{Q_{1,m}(w)}
&=&\frac{1+p_{1,m,1}w}
{1+q_{1,m,1}w+\dots+q_{1,m,m}w^m}\\
&=&\zetatilde(w)+O(w^{m+2})
\end{eqnarray}
It can be checked that
\begin{equation}
p_{1,m,1}=-\frac{\tau_{m+1}}{\tau_{m}}
\end{equation}
and hence, according to \eqref{RH1ver2}, for $m\rightarrow \infty$
\begin{equation}
p_{1,m,1}\rightarrow \frac{3}{2}
\label{limp1ml}\end{equation}
and, respectively,
\begin{equation}
P_{1,m}(w)\rightarrow 1+\tfrac{3}{2}w=1-\frac{w}{\w_1}.
\label{limP1}\end{equation}

Now let us consider the general case. Let
let $P_{l,m}(w)$ and $Q_{l,m}(w)$ be such polynomials that
\begin{eqnarray}
\zetatilde(w)\approx\frac{P_{l,m}(w)}{Q_{l,m}(w)}
&=&\frac{1+p_{l,m,1}w+\dots +p_{l,m,l}w^l}
{1+q_{l,m,1}w+\dots+q_{l,m,m}w^m}\\
&=&\zetatilde(w)+O(w^{l+m+1})
\end{eqnarray}
According to a theorem of de Montessue \cite{Montessue1902,Montessue1905}
(see also  \cite{BakerPade}) for every $l$, subhypothesis $\RH_l$ implies the following
generalization of \eqref{limP1}:
\emph{
 For all $l$ for $m\rightarrow \infty$
\begin{equation}
P_{l,m}(w)\rightarrow \prod_{k=1}^l\left(1-\frac{w}{\w_l}\right).
\label{limPlm}\end{equation}
}

We are going to deal only with the leading coefficient of $P_{l,m}(w)$
for which \eqref{limPlm} implies
 the following generalization of \eqref{limp1ml}:

\

{\bf Subhypothesis $\mathbf{RH_\emph{l}^w}$ (version 1)}.  \emph{For $m\rightarrow \infty$
\begin{eqnarray}
p_{l,m,l}\rightarrow \W_l
\label{limplml}\end{eqnarray}
where
\begin{eqnarray}
\W_l=\prod_{k=1}^l\left(-\frac{1}{\w_l}\right)=
\prod_{k=1}^l\frac{2k+1}{2k}.
\label{defW}\end{eqnarray}
}

\

\subsection{Back to the Riemann Hypothesis}

Each subhypothesis $\RH^\mathrm{w}_l$ is, formally, weaker than
the corresponding subhypothesis
$\RH_l$, nevertheless, taken together the subhypotheses $\RH^\mathrm{w}_l$ are equivalent
to the subhypotheses
$\RH_l$, and thus
we have

\

{\bf Riemann Hypothesis (version 5).} \emph{For every $m$
the subhypothesis $\RH^\mathrm{w}_l$ is true.}

\

In order to see why it is so, suppose that the Riemann Hypothesis
isn't valid, and let $\check{z}$ be a non-trivial zero of $\zeta(z)$
with $\Re(\check z)<\frac{1}{2}$.
Then $\check w=\frac{\check z}{1-\check z}$ is a non-trivial zero of $\zetatilde(w)$
with $\abs{\check w}<1$. In the closed circle $\abs{w}\le \abs{\check w}$
there are only finitely many zeros of $\zetatilde(w)$; let us denote them by
$\check w_1$, \dots,  $\check w_l$. By the above cited  theorem of de Montessue,
for $m\rightarrow \infty$
\begin{equation}
P_{l,m}(w)\rightarrow \prod_{k=1}^l\left(1-\frac{w}{\check w_l}\right)
\end{equation}
and hence
\begin{eqnarray}
p_{l,m,l}\rightarrow \prod_{k=1}^l\left(-\frac{1}{\check w_l}\right).
\end{eqnarray}
It is easy to see that
\begin{eqnarray}
\left|\prod_{k=1}^l\left(-\frac{1}{\check w_l}\right)\right|>
\left|\prod_{k=1}^l\left(-\frac{1}{\w_l}\right)\right|=\abs{\W_l}
\end{eqnarray}
which gives the required contradiction with \eqref{limplml}.

\subsection{More Determinants}

An explicit expression for $p_{l,m,l}$ can be
given (Jacobi \cite{Jacobi1846}, see also \cite{BakerPade}):
\begin{equation}
p_{l,m,l}=\frac{\det({L_{l,m+1}}(\zetastar))}{\det({L_{l,m}}(\zetastar))}
\label{plml}\end{equation}
where
\begin{equation}
L_{l,m}(f)=
\begin{pmatrix}
\theta_{l}(f)&\theta_{l-1}(f)&\dots&\theta_{l-m+1}(f)\\
\theta_{l+1}(f)&\theta_{l}(f)&\dots&\theta_{l-m+2}(f)\\
\vdots&\vdots&\ddots&\vdots\\
\theta_{l+m-1}(f)&\theta_{l+m-2}(f)&\dots&\theta_{l}(f)
\end{pmatrix}
\end{equation}
with $\theta_0(f)=1$ and $\theta_j(f)=0$   for  $j<0$.

In terms of these matrices we have the following counterpart of
\eqref{det1R1}:

\

{\bf  $\mathbf{RH_\emph{l}^w}$ (version 2)}.  \emph{For $m\rightarrow \infty$
\begin{equation}
\det(L _{l,m}(\zetastar))\,
=\,\W_l^m(R_l(\zetastar)+o(1)).
\label{detlR1}\end{equation}
with some constant $R_l(\zetastar)$}.

\

In order to pass from \eqref{RH1ver2}, the second version of $\RH_1$, to
\eqref{RH1ver3}, the third version of $\RH_1$, we needed the inequality
\eqref{R1pos} which we got from the numerical value \eqref{R1num}.
However, it is easy to see that $\RH_1(f)$ implies the inequality $R_1(f)>0$
for every function $f$ satisfying condition \eqref{f01}. Namely, by analogy with
\eqref{R1tilde} $\RH_1(f)=-\frac{1}{\z_1f'(\z_1)}$ and $f'(\z_1)>0$
because $\z_1=-\frac{2}{3}$ is the least (in absolute value) zero of $f(z)$. By a similar
argument, for every $l$   the inequality $R_l(f)>0$ is implied by $\RH_l(f)$,
and we have

\

{\bf $\mathbf{RH_\emph{l}^w}$ (version 3)}.  \emph{
\begin{equation}
\lim_{m\rightarrow \infty}\left(\det(L _{l,m}(\zetastar))\right)^\frac{1}{m}\,
=\,\W_l.
\label{detlRl}\end{equation}
}

\subsection{More Eigenvalues}

By analogy with \eqref{lambdas}, we have the representation
\begin{equation}
\det(L _{l,m}(f))=\lambda_{l,m,1}(f)\lambda_{l,m,2}(f)\dots
\lambda_{l,m,m}(f)
\label{llambdas}\end{equation}
where $\lambda_{l,m,1}(f), \ \lambda_{l,m,2}(f),\ \dots,
 \ \lambda_{l,m,m}(f)$
are the eigenvalues of the matrix $L_{l,m}(f)$.
Then, we have

\

{\bf $\mathbf{RH_\emph{l}^w}$ (version 4).}
\begin{equation}
\lim_{m\rightarrow \infty}
\left(\prod_{n=1}^m{\lambda_{l,m,n}(\zetastar)}\right)^\frac{1}{m}=
\W_l.
\end{equation}

\

The (multi)set $\{\lambda_{l,m,1}(f),\ \lambda_{l,m,2}(f),\dots
\lambda_{l,m,m}(f)\}$ will be called
the $\lambda$-\emph{spectrum} of the function $f$ and will be
denoted~$\lS_{l,m}(f)$.

\subsection{More about positions of eigenvalues}

\begin{figure}[p]
\includegraphics[width=.47\textwidth]{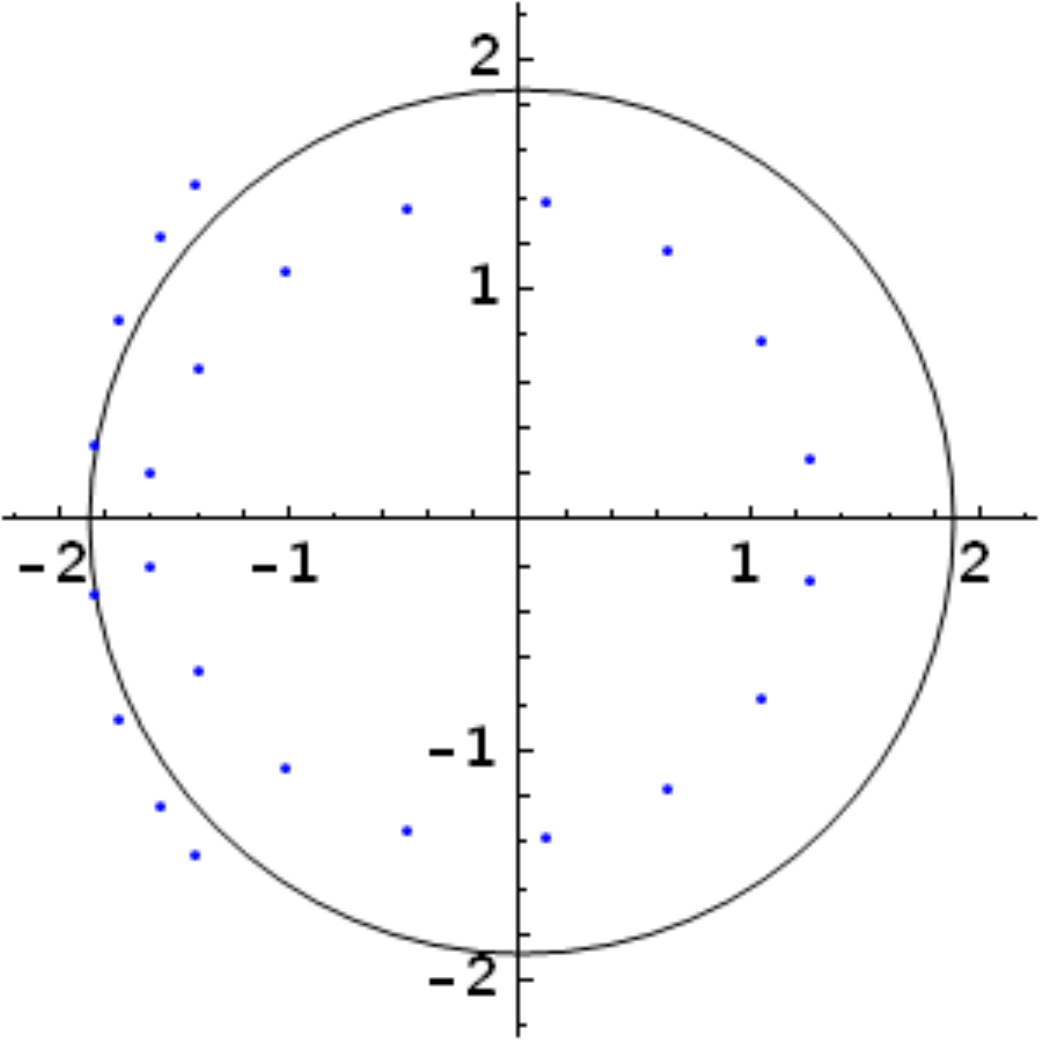}\hfill
\includegraphics[width=.47\textwidth]{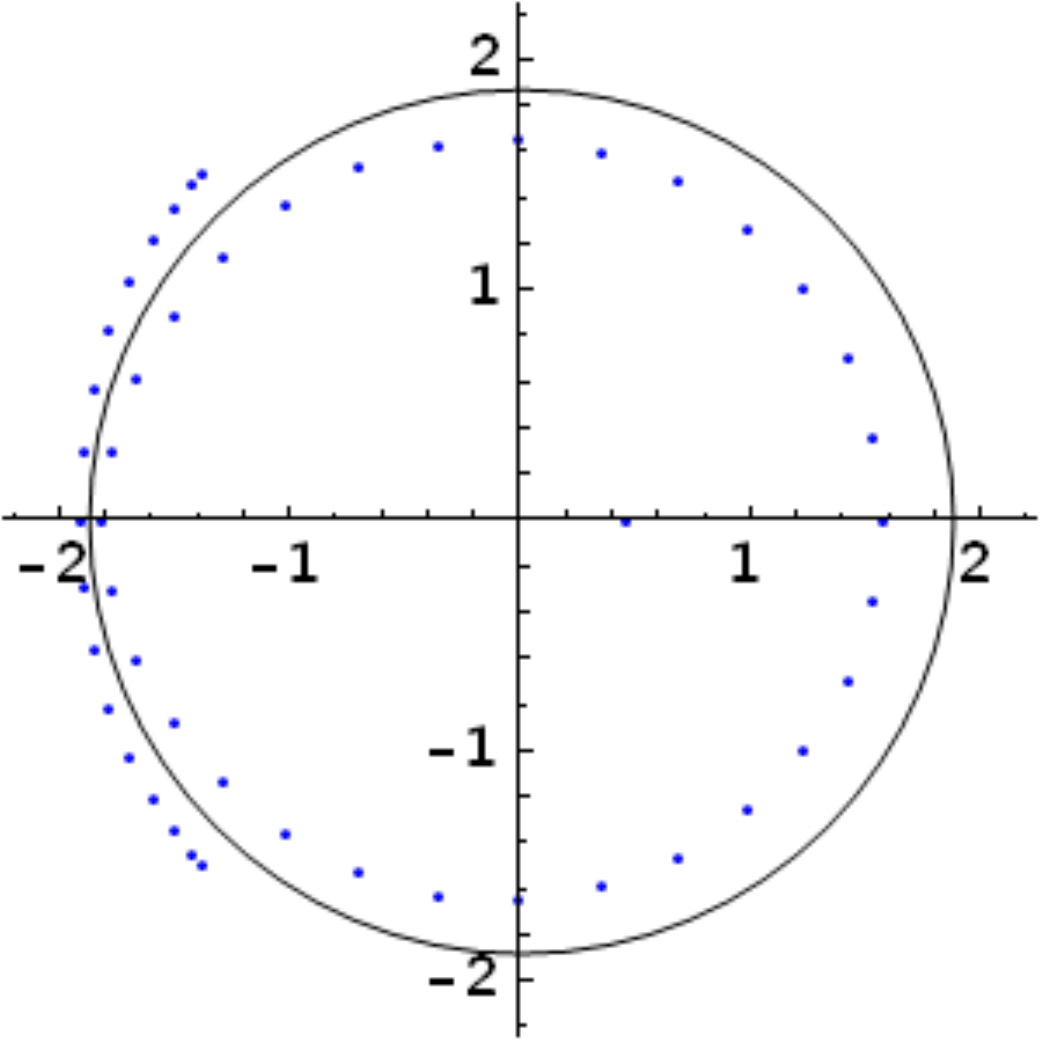}\\
\parbox{.47\textwidth}{\caption{$\lS_{2,24}(\zetastar)$}\label{pic224}}
\hfill\parbox{.47\textwidth}{\caption{$\lS_{2,48}(\zetastar)$}\label{pic248}}
\end{figure}
\begin{figure}[p]
\includegraphics[width=.47\textwidth]{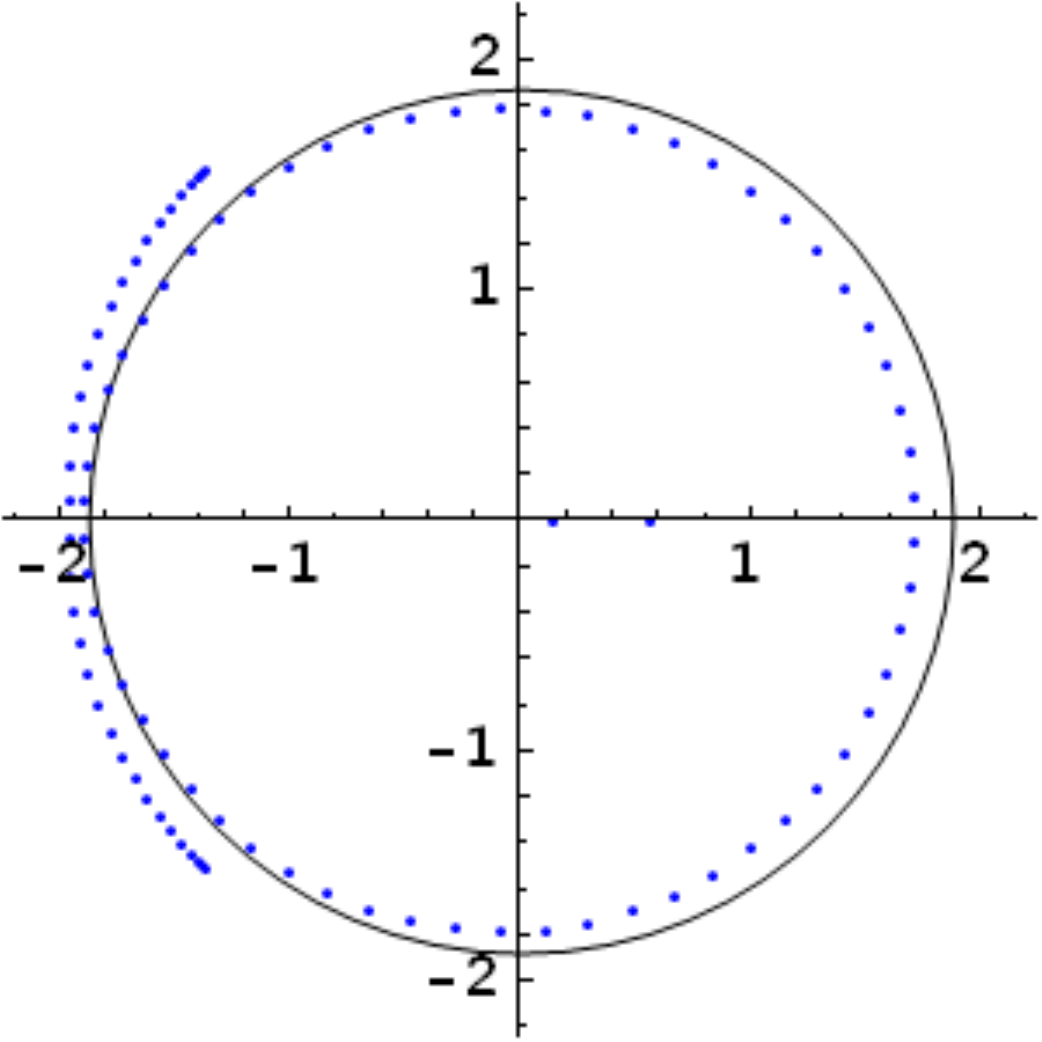}\hfill
\includegraphics[width=.47\textwidth]{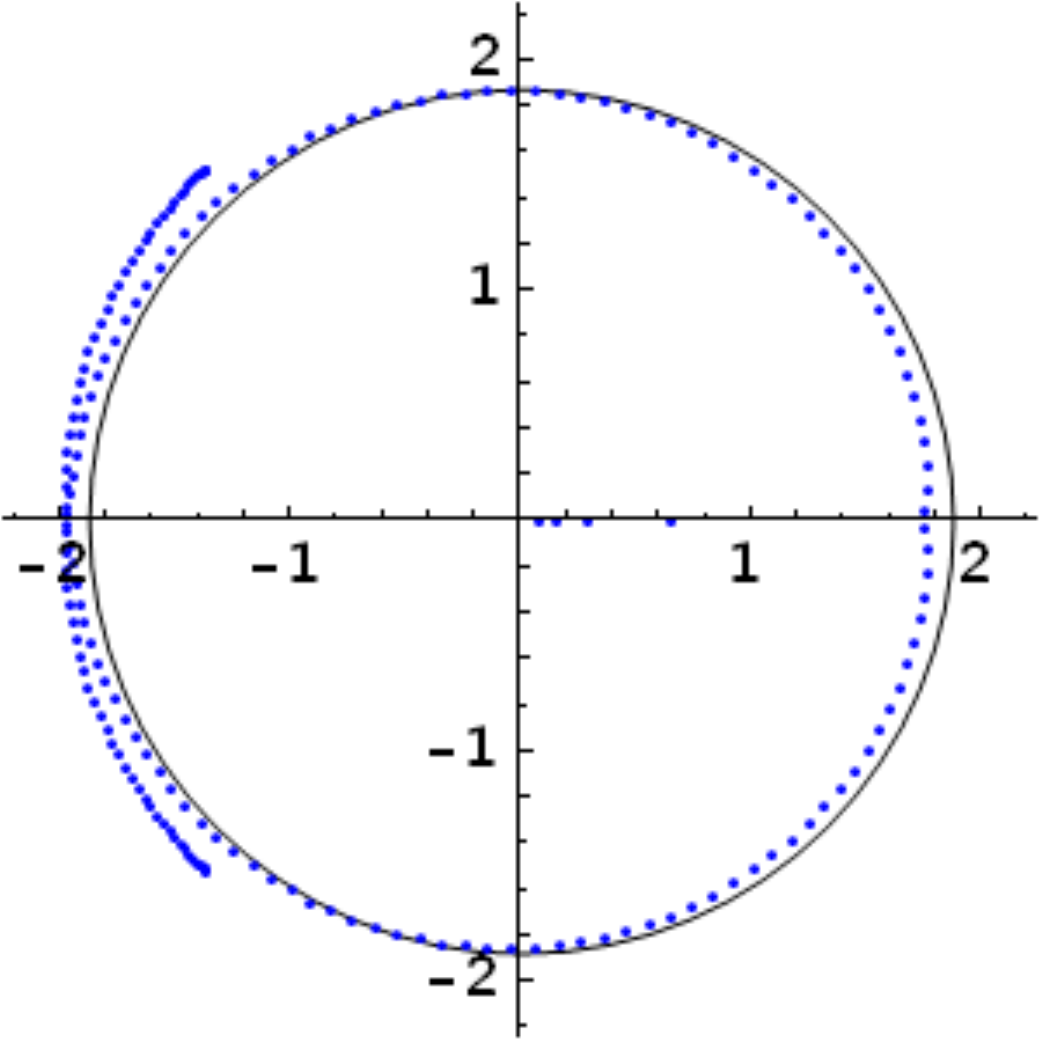}\\
\parbox{.47\textwidth}{\caption{$\lS_{2,96}(\zetastar)$}\label{pic296}}\hfill
\parbox{.47\textwidth}{\caption{$\lS_{2,192}(\zetastar)$}\label{pic2192}}
\end{figure}
Figures \ref{pic224}--\ref{pic2192}  show  spectra $\lS_{2,m}(\zetastar)$
for $m=24,\ 48,\ 96, \ 192$
  respectively.
%Figure \ref{pic21192} shows  the union $\cup_{m=1}^{192}\lS_{2,m}(\zetastar)$.
An
animation showing the  $\lS_{2,1}(\zetastar)$, $\lS_{2,2}(\zetastar)$, \dots in succession can
be downloaded from~\cite{hidden}.

We see that $\lS_{2,m}(\zetastar)$ consists of the arrow, the bow (now looking
into the
opposite direction), and a new element, looking like a circle,
which will be called \emph{orbit}.
The animation shows that the orbit has,
on its right-hand side, a
\emph{rendezvous}
with the arrow  and,  on its left-hand side,
another \emph{rendezvous}
with the bow.

\begin{figure}[p]
\includegraphics[width=.47\textwidth]{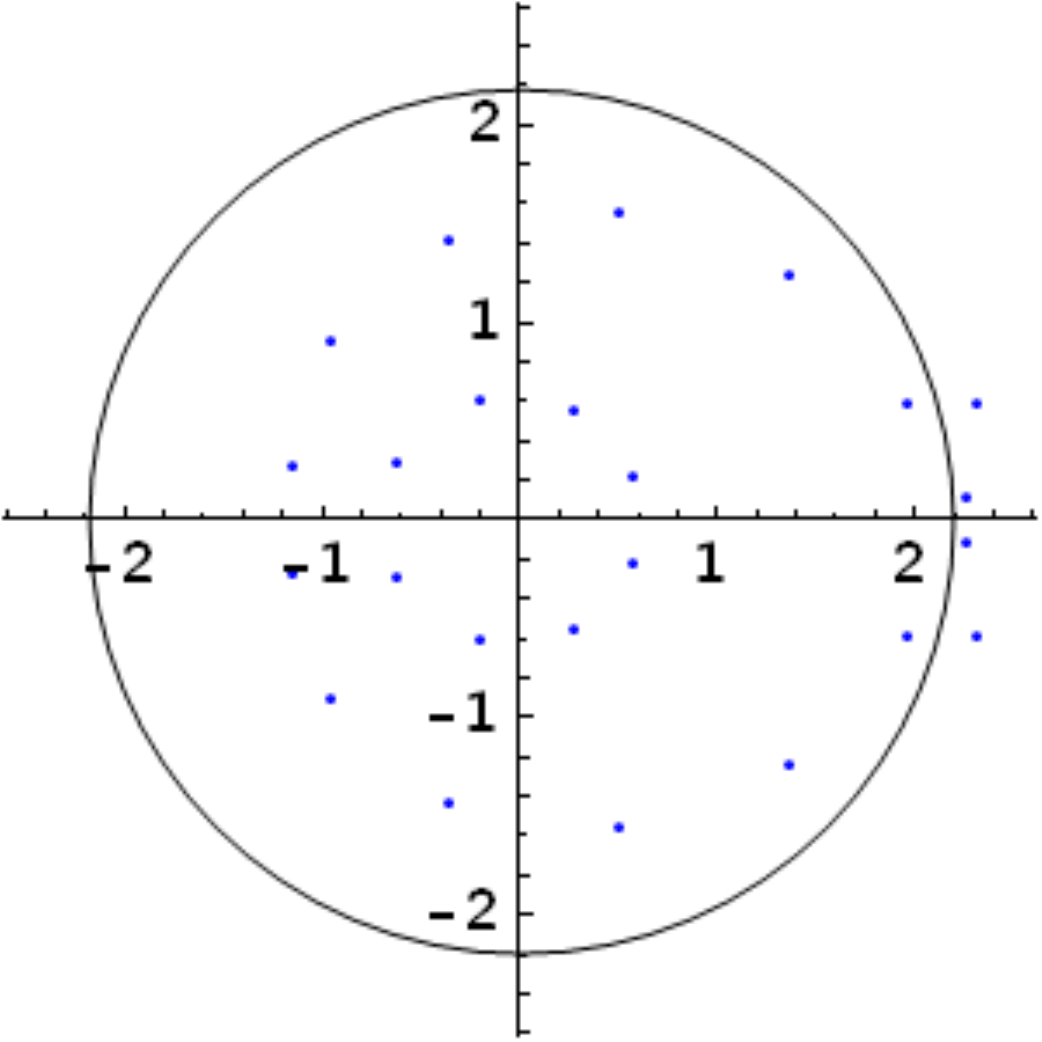}\hfill
\includegraphics[width=.47\textwidth]{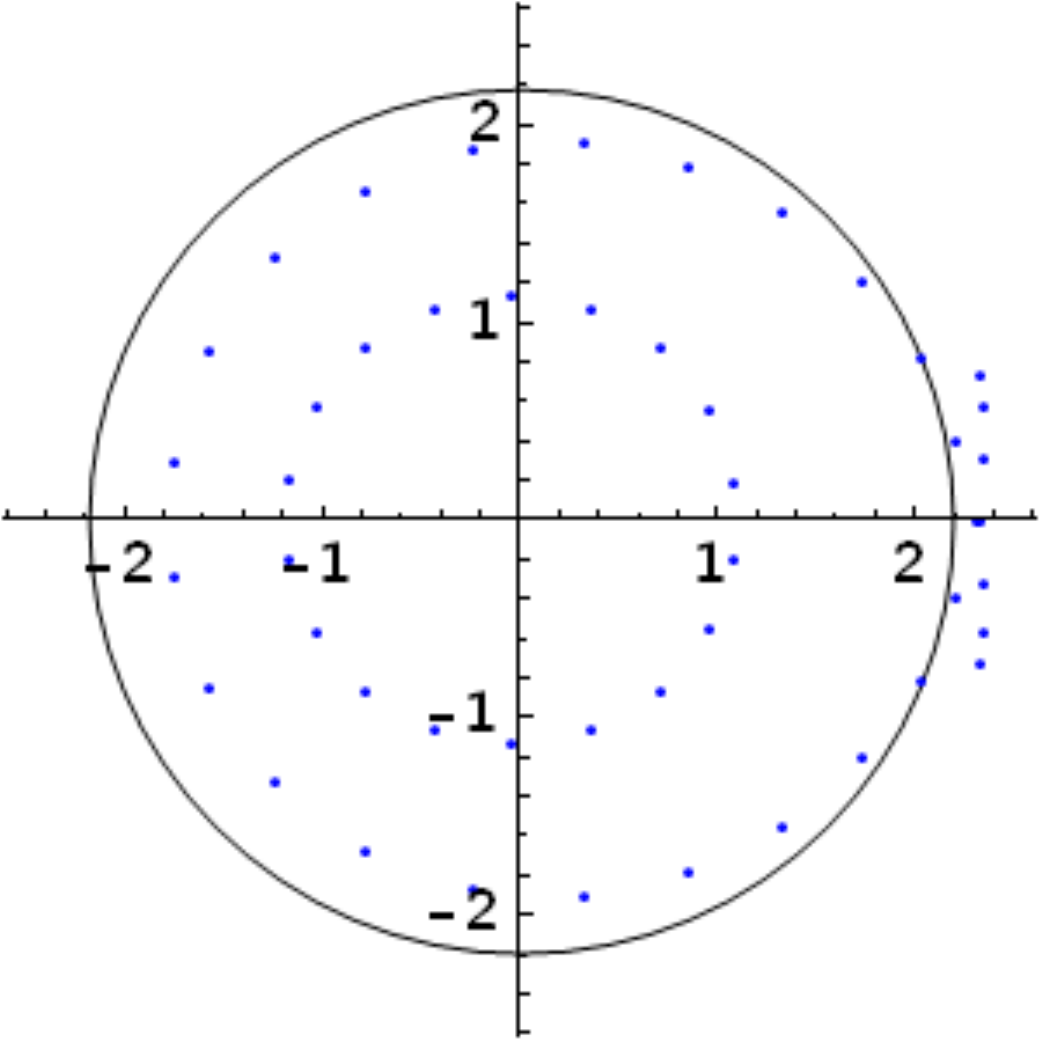}\\
\parbox{.47\textwidth}{\caption{$\lS_{3,24}(\zetastar)$}\label{pic324}}
\hfill\parbox{.47\textwidth}{\caption{$\lS_{3,48}(\zetastar)$}\label{pic348}}
\end{figure}
\begin{figure}[p]
\includegraphics[width=.47\textwidth]{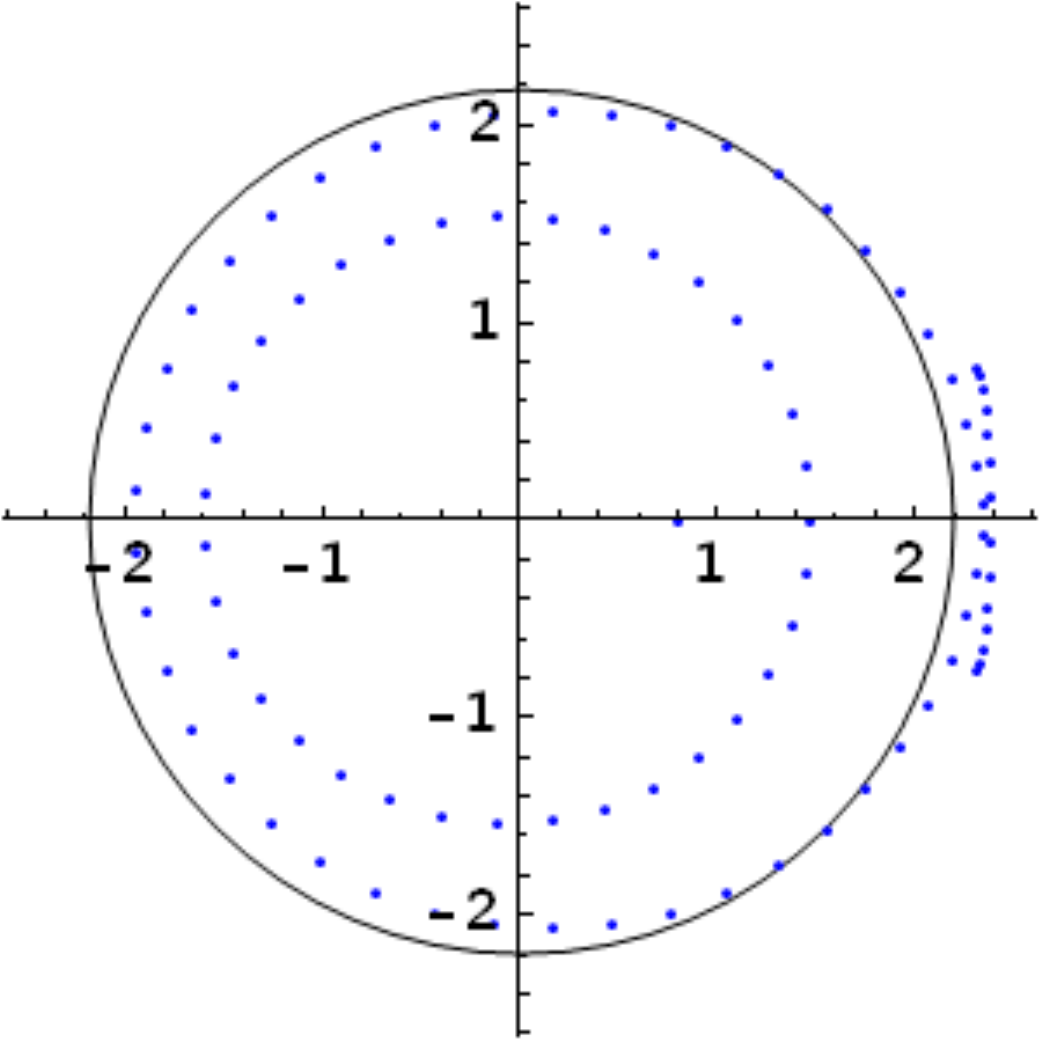}\hfill
\includegraphics[width=.47\textwidth]{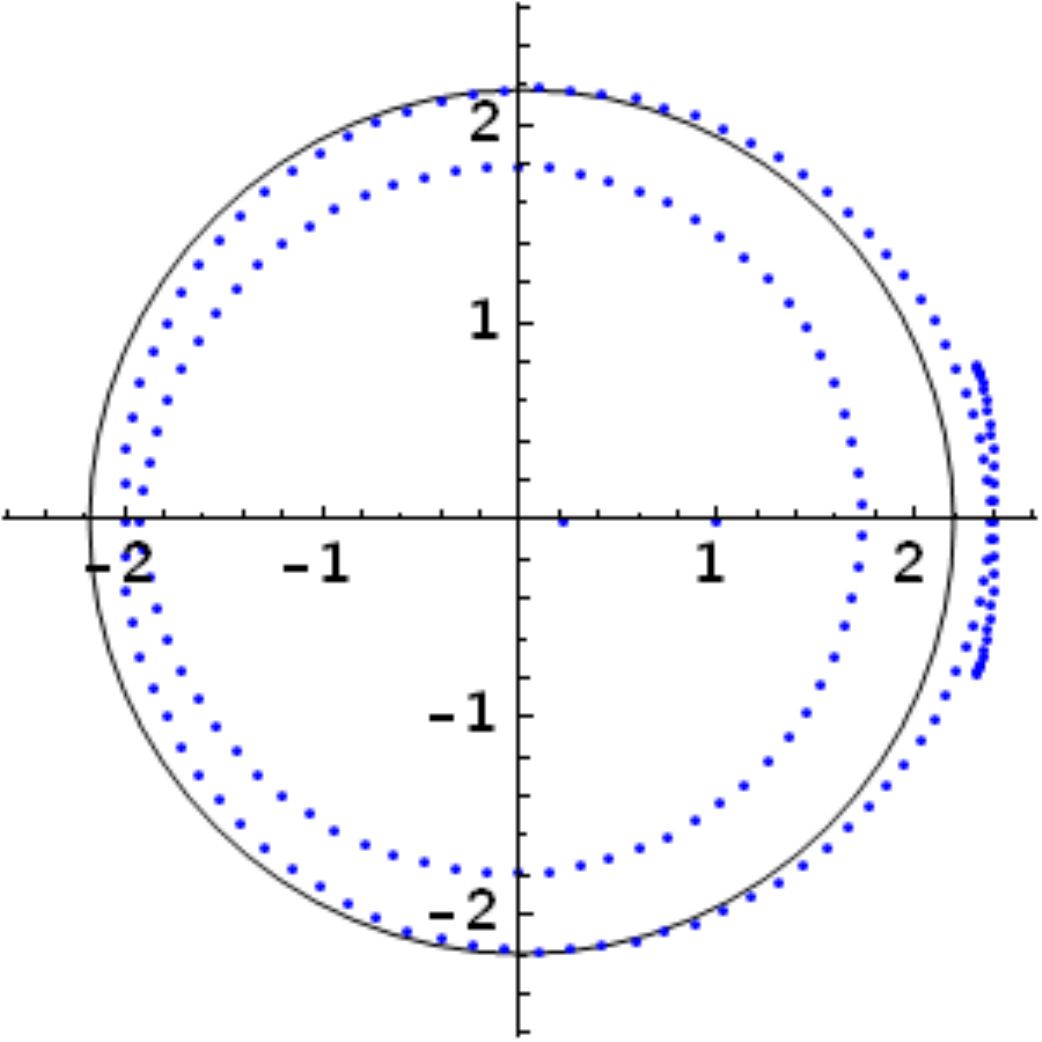}\\
\parbox{.47\textwidth}{\caption{$\lS_{3,96}(\zetastar)$}\label{pic396}}\hfill
\parbox{.47\textwidth}{\caption{$\lS_{3,192}(\zetastar)$}\label{pic3192}}
\end{figure}
Figures \ref{pic324}--\ref{pic3192}  show  spectra $\lS_{3,m}(\zetastar)$
for $m=24,\ 48,\ 96, \ 192$
  respectively.
%Figure \ref{pic3192} shows  the union $\cup_{m=1}^{192}\lS_{3,m}(\zetastar)$.
An animation showing the  $\lS_{3,1}(\zetastar)$, $\lS_{3,2}(\zetastar)$, \dots in succession can
be downloaded from \cite{hidden}.

We see that $\lS_{3,m}(\zetastar)$ consists of the arrow, the bow (now rather rudimentary
and looking into the same direction as in the case of $\lS_{1,m}$), and
two orbits which constitute \emph{target}.
The animation shows that  the inner orbit has, on its right-hand side,
the rendezvous
with the arrow  and has, on its left-hand side, the rendezvous
with the outer orbit. In its turn, the outer orbit has, on its left-hand side,
the rendezvous with the inner orbit
  and has, on its right-hand side, the rendezvous with the bow.

\begin{figure}[p]
\includegraphics[width=.47\textwidth]{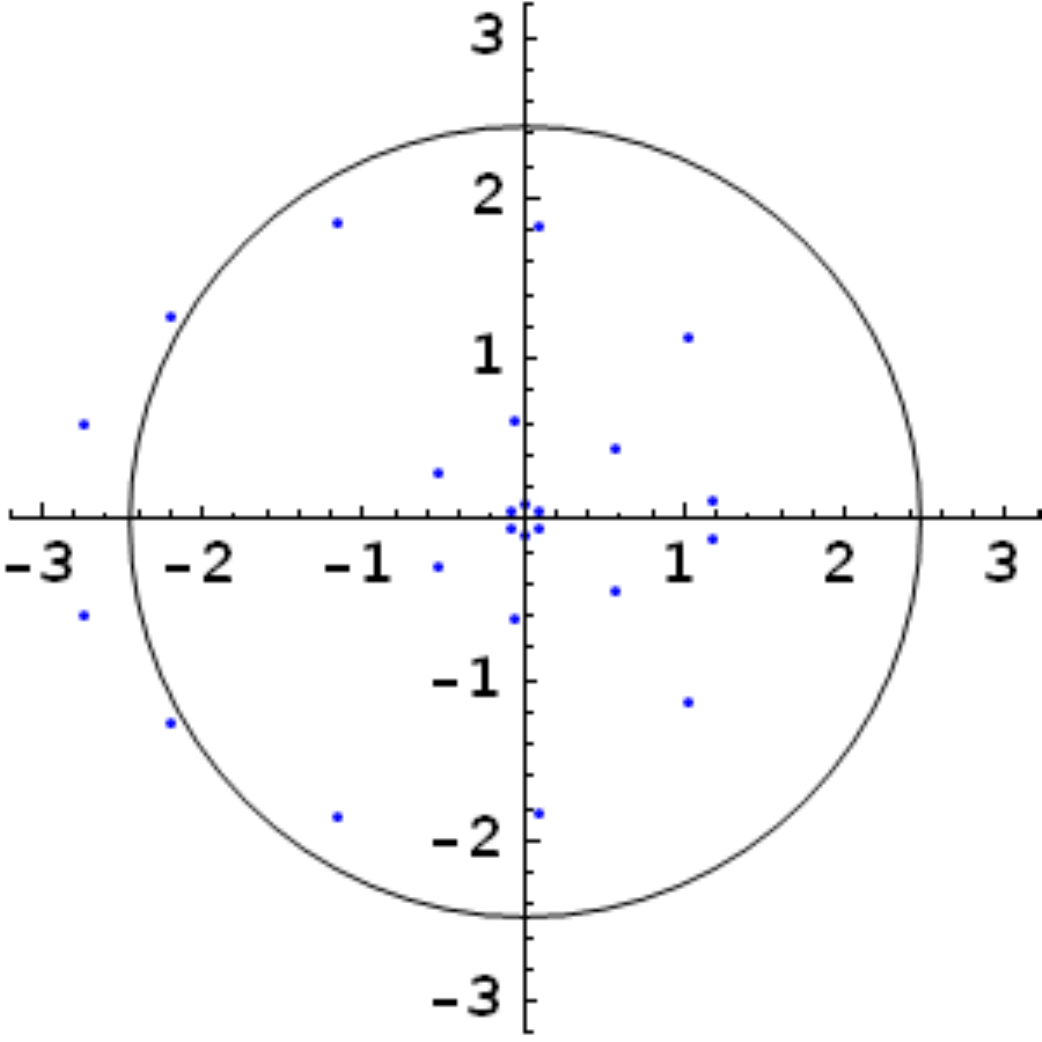}\hfill
\includegraphics[width=.47\textwidth]{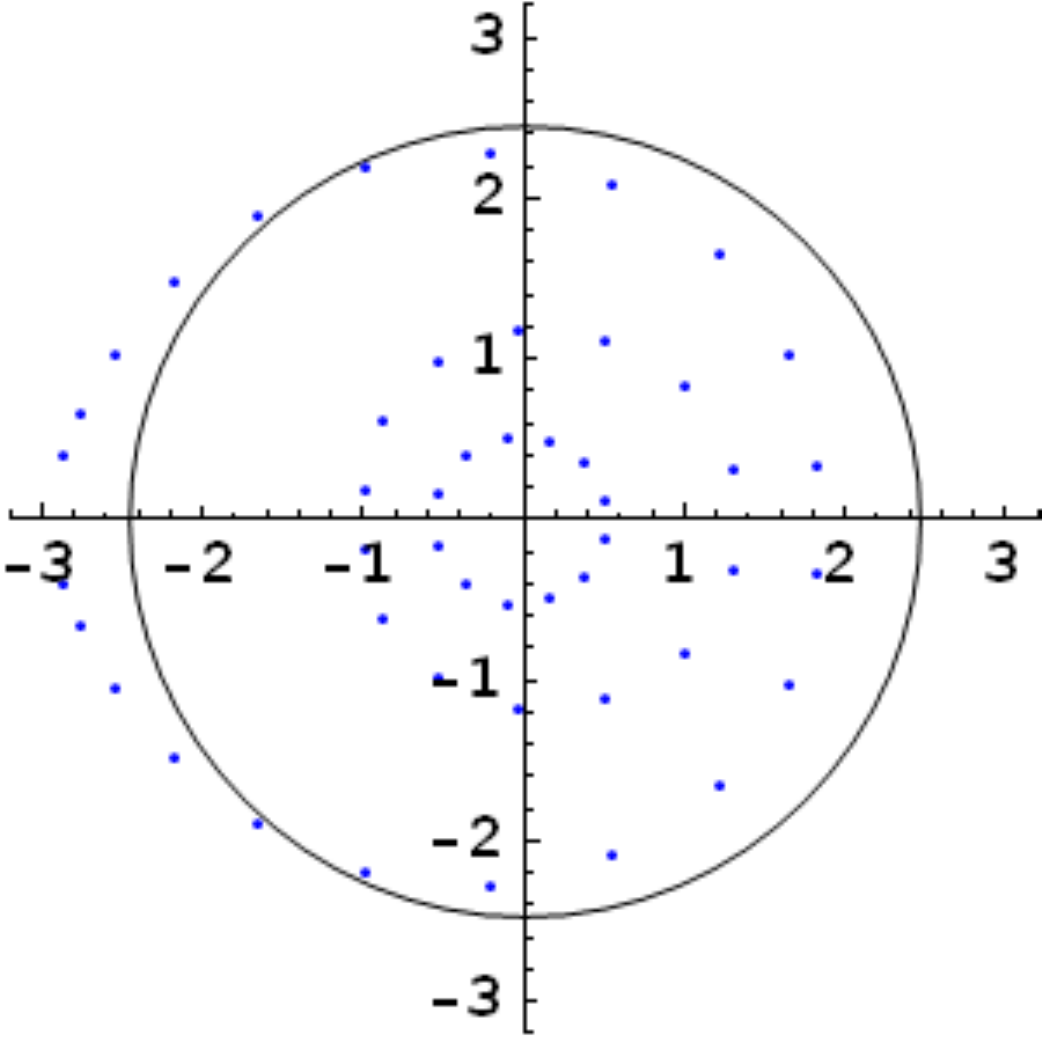}\\
\parbox{.47\textwidth}{\caption{$\lS_{4,24}(\zetastar)$}\label{pic424}}
\hfill\parbox{.47\textwidth}{\caption{$\lS_{4,48}(\zetastar)$}\label{pic448}}
\end{figure}
\begin{figure}[p]
\includegraphics[width=.47\textwidth]{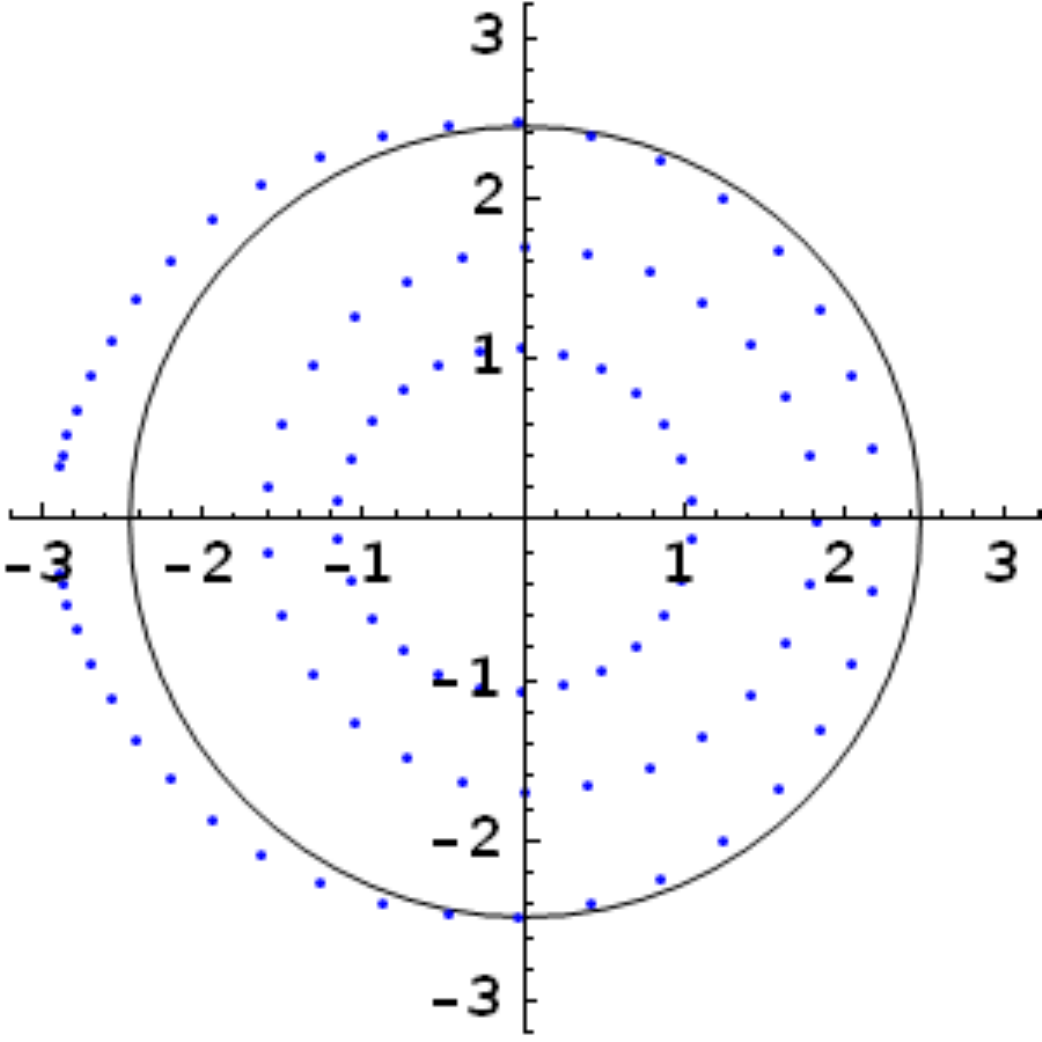}\hfill
\includegraphics[width=.47\textwidth]{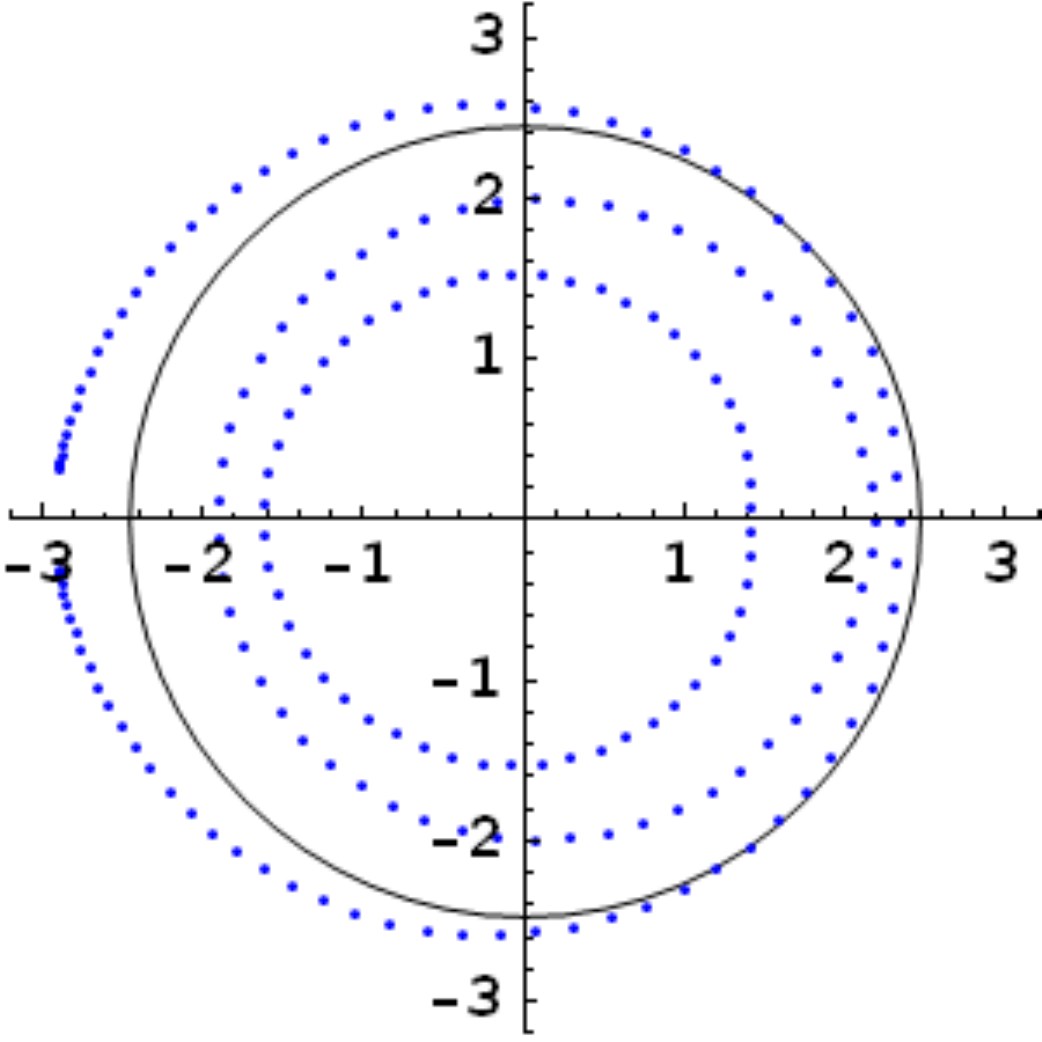}\\
\parbox{.47\textwidth}{\caption{$\lS_{4,96}(\zetastar)$}\label{pic496}}\hfill
\parbox{.47\textwidth}{\caption{$\lS_{4,192}(\zetastar)$}\label{pic4192}}
\end{figure}
 One might expect that the target of the spectra $\lS_{4,m}(\zetastar)$ would consist of
 three orbits but this is not the case.
 Figures \ref{pic424}--\ref{pic4192}   show  spectra $\lS_{4,m}(\zetastar)$
for $m=24,\ 48,\ 96, \ 192$
  respectively.
%Figure \ref{pic41192} shows  the union $\cup_{m=1}^{192}\lS_{4,m}(\zetastar)$.
An animation showing the  $\lS_{4,1}(\zetastar)$, $\lS_{4,2}(\zetastar)$, \dots in succession can
be downloaded from \cite{hidden}.
We see that the entire structure of the spectra $\lS_{4,m}(\zetastar)$ is the same
as in the case of the spectra $\lS_{3,m}(\zetastar)$ but the bow isn't rudimentary any longer,
on the contrary, it is almost a circle.

\begin{figure}[p]
\includegraphics[width=.47\textwidth]{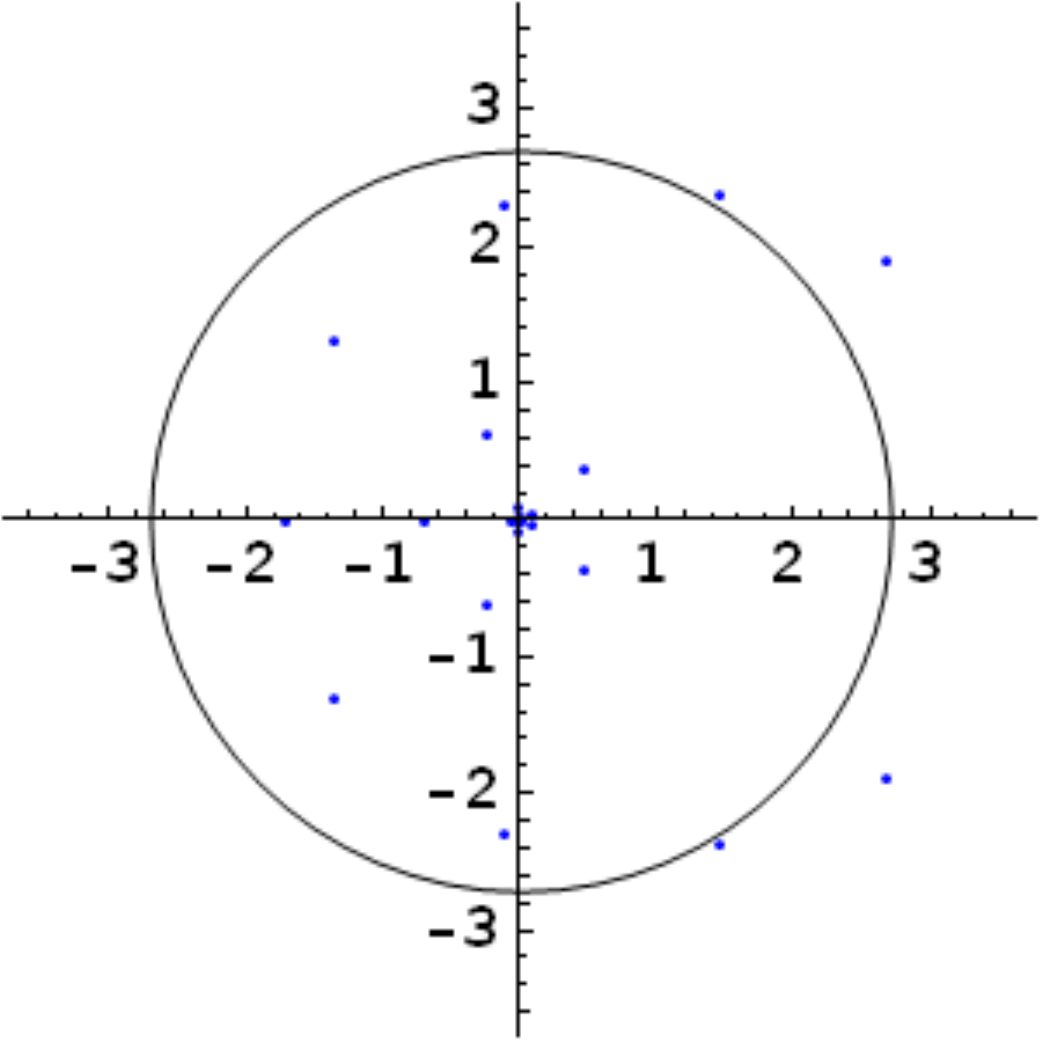}\hfill
\includegraphics[width=.47\textwidth]{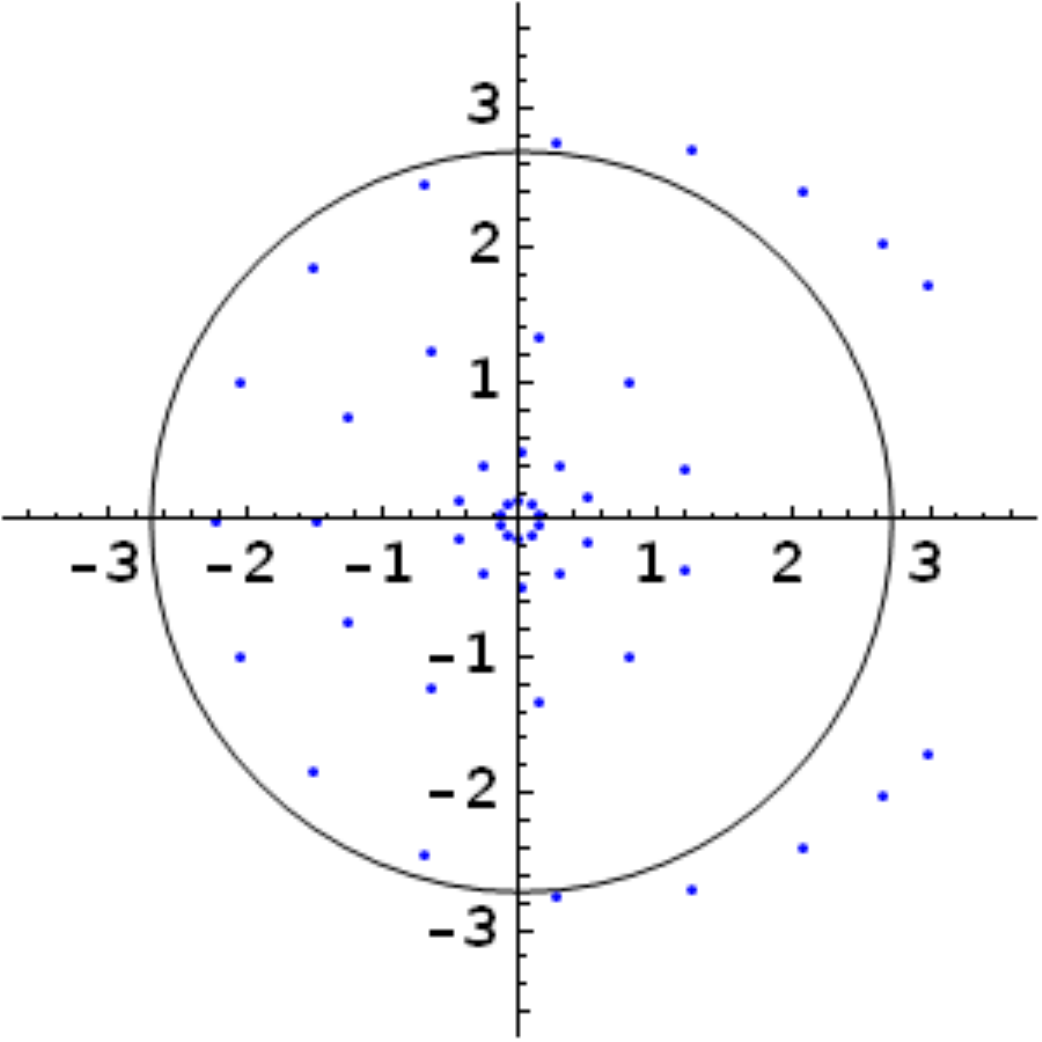}\\
\parbox{.47\textwidth}{\caption{$\lS_{5,24}(\zetastar)$}\label{pic524}}
\hfill\parbox{.47\textwidth}{\caption{$\lS_{5,48}(\zetastar)$}\label{pic548}}
\end{figure}
\begin{figure}[p]
\includegraphics[width=.47\textwidth]{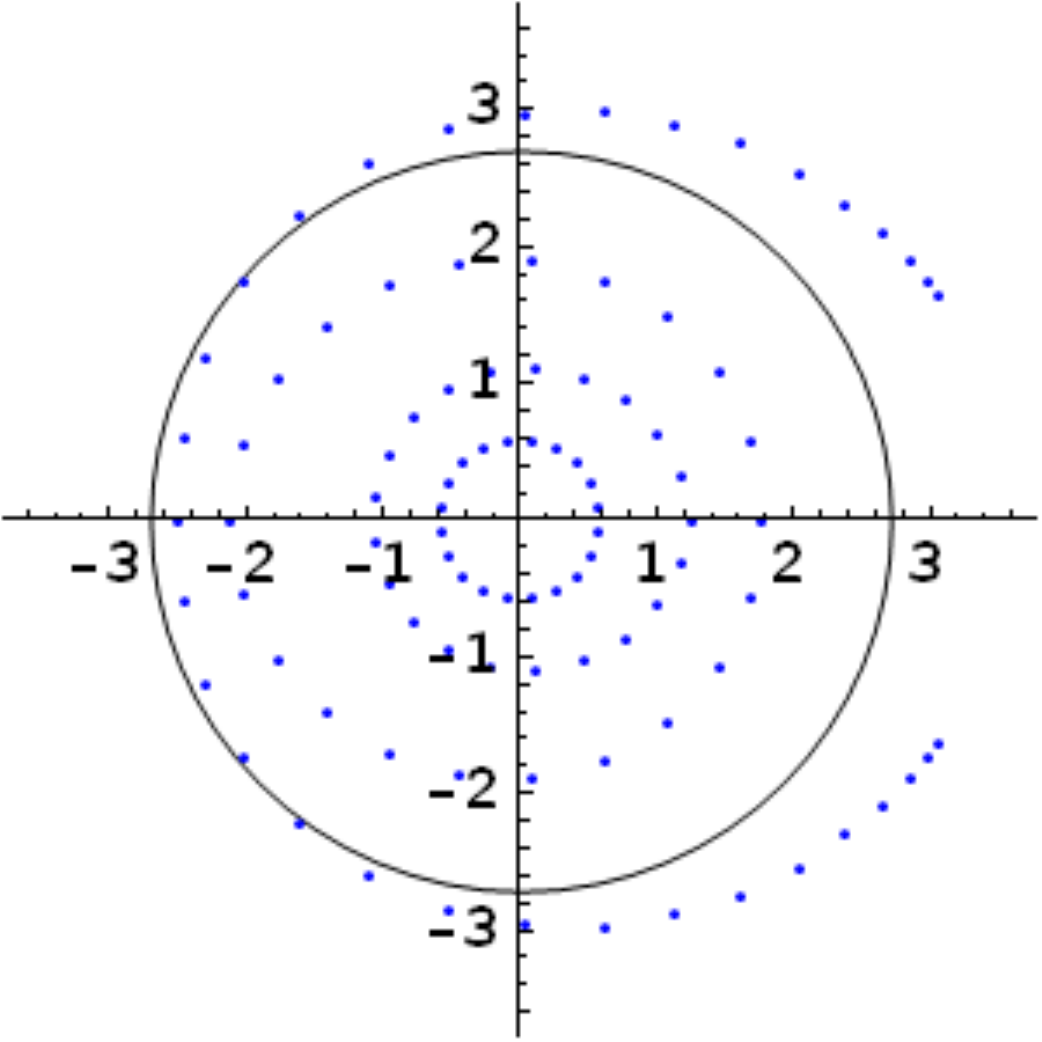}\hfill
\includegraphics[width=.47\textwidth]{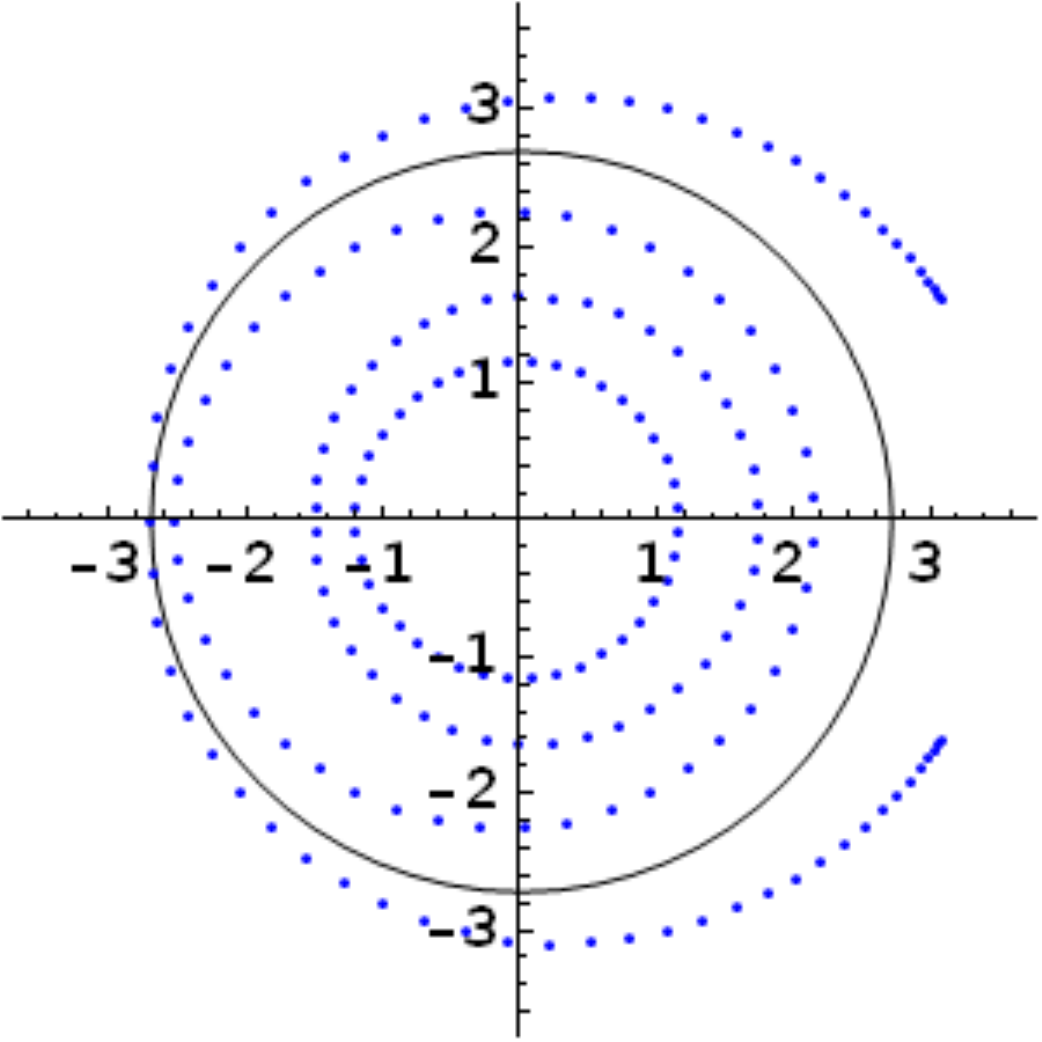}\\
\parbox{.47\textwidth}{\caption{$\lS_{5,96}(\zetastar)$}\label{pic596}}\hfill
\parbox{.47\textwidth}{\caption{$\lS_{5,192}(\zetastar)$}\label{pic5192}}
\end{figure}
The third orbit in the target appears in the spectra $\lS_{5,m}(\zetastar)$.
 Figures \ref{pic524}--\ref{pic5192} show     spectra $\lS_{5,m}(\zetastar)$
for $m=24,\ 48,\ 96, \ 192$
  respectively.
%Figure \ref{pic51192} shows  the union $\cup_{m=1}^{192}\lS_{5,m}(\zetastar)$.
An animation showing the  $\lS_{5,1}(\zetastar)$, $\lS_{5,2}(\zetastar)$, \dots in succession can
be downloaded from \cite{hidden}.

\begin{figure}[p]
\includegraphics[width=.47\textwidth]{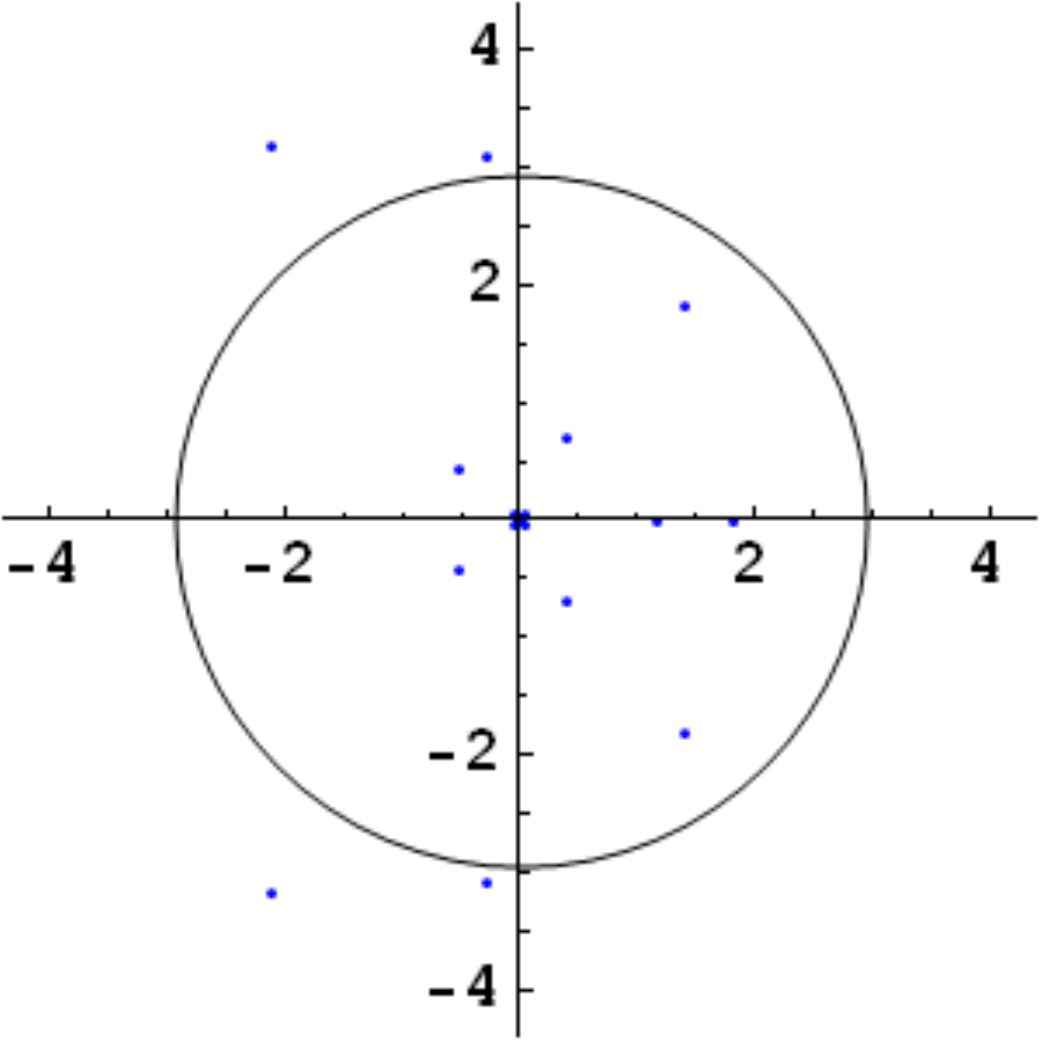}\hfill
\includegraphics[width=.47\textwidth]{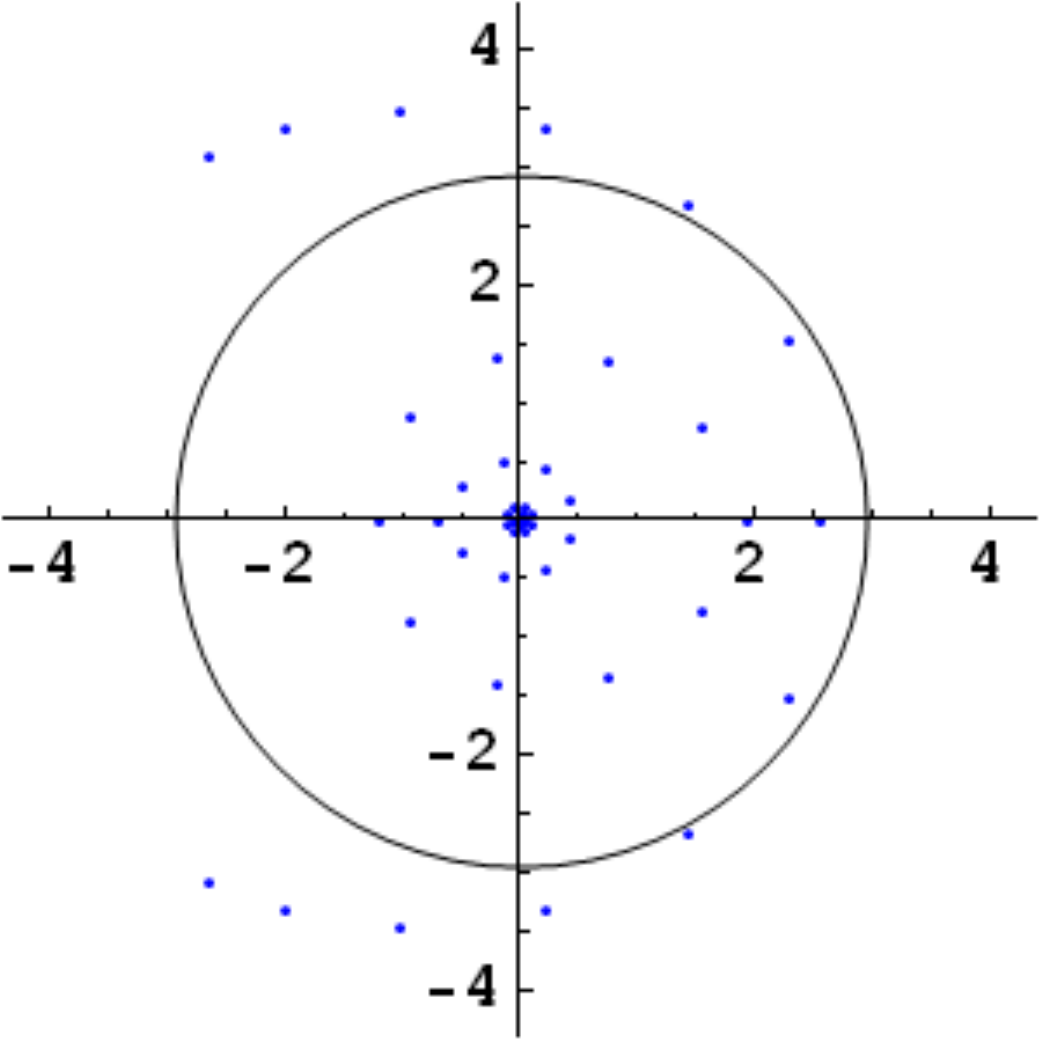}\\
\parbox{.47\textwidth}{\caption{$\lS_{6,24}(\zetastar)$}\label{pic624}}
\hfill\parbox{.47\textwidth}{\caption{$\lS_{6,48}(\zetastar)$}\label{pic648}}
\end{figure}
\begin{figure}[p]
\includegraphics[width=.47\textwidth]{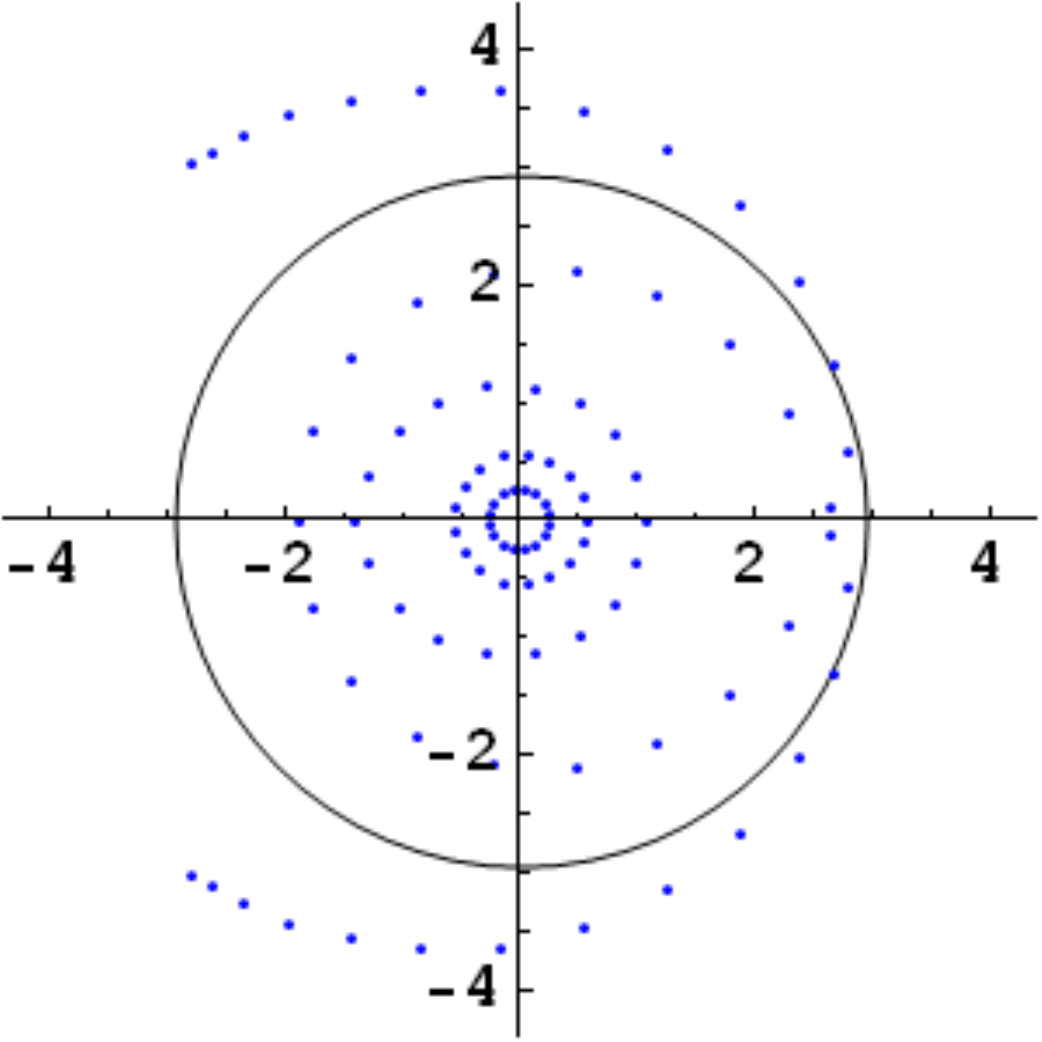}\hfill
\includegraphics[width=.47\textwidth]{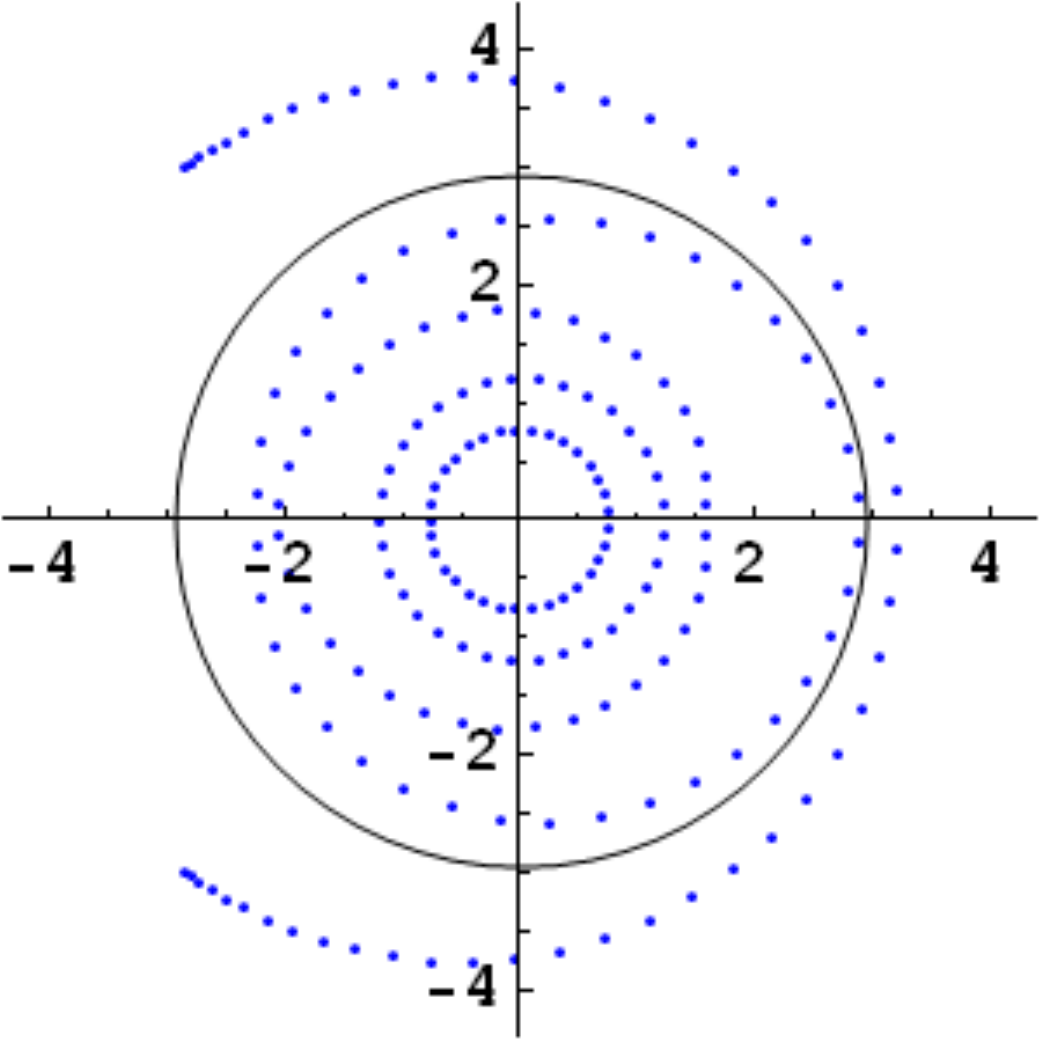}\\
\parbox{.47\textwidth}{\caption{$\lS_{6,96}(\zetastar)$}\label{pic696}}\hfill
\parbox{.47\textwidth}{\caption{$\lS_{6,192}(\zetastar)$}\label{pic6192}}
\end{figure}
Similar, the fourth orbit in the target appears in the spectra $\lS_{6,m}(\zetastar)$.
 Figures \ref{pic624}--\ref{pic6192}  show  spectra $\lS_{6,m}(\zetastar)$
for $m=24,\ 48,\ 96, \ 192$
  respectively.
%Figure \ref{6192} shows  the union $\cup_{m=1}^{192}\lS_{6,m}(\zetastar)$.
An animation showing the  $\lS_{6,1}(\zetastar)$, $\lS_{6,2}(\zetastar)$, \dots in succession can
be downloaded from \cite{hidden}.

\begin{figure}[p]
\includegraphics[width=.47\textwidth]{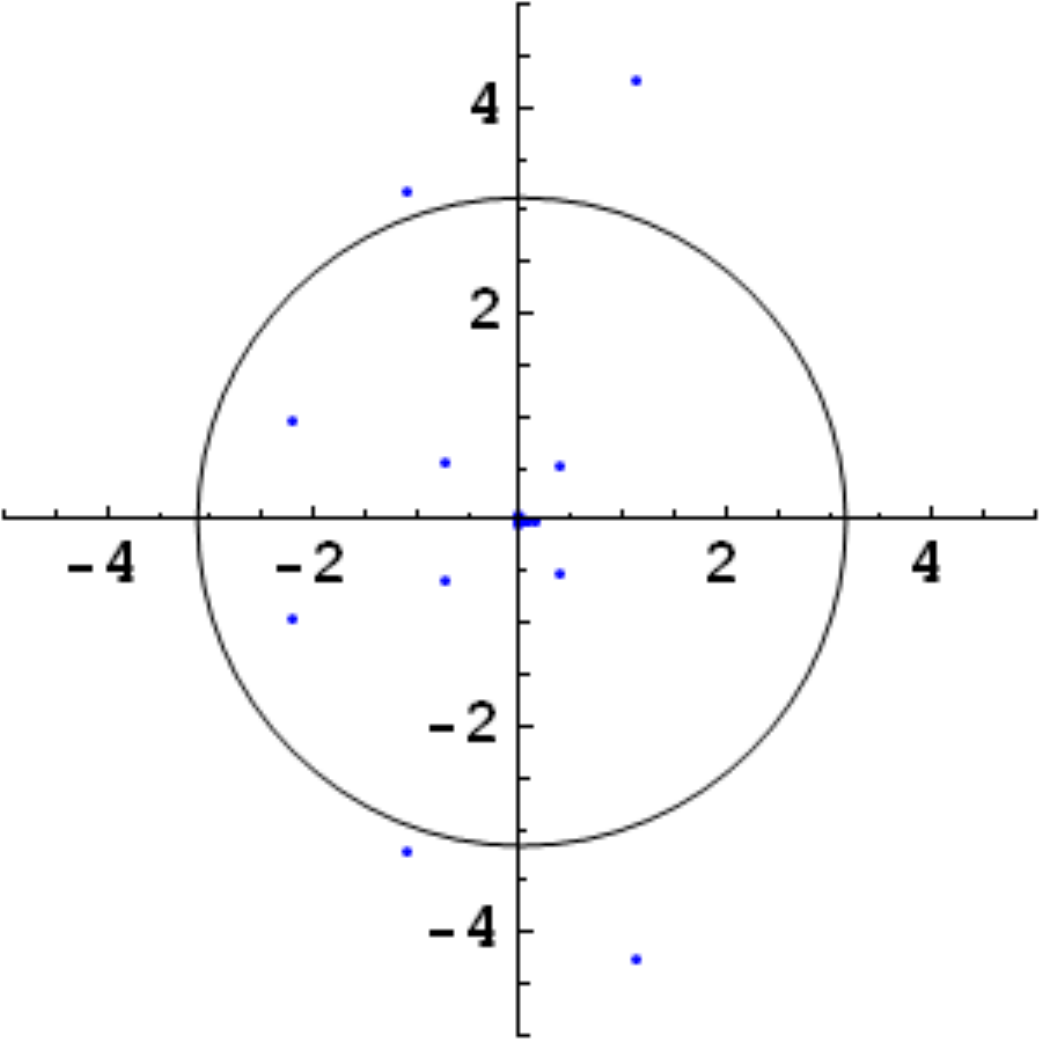}\hfill
\includegraphics[width=.47\textwidth]{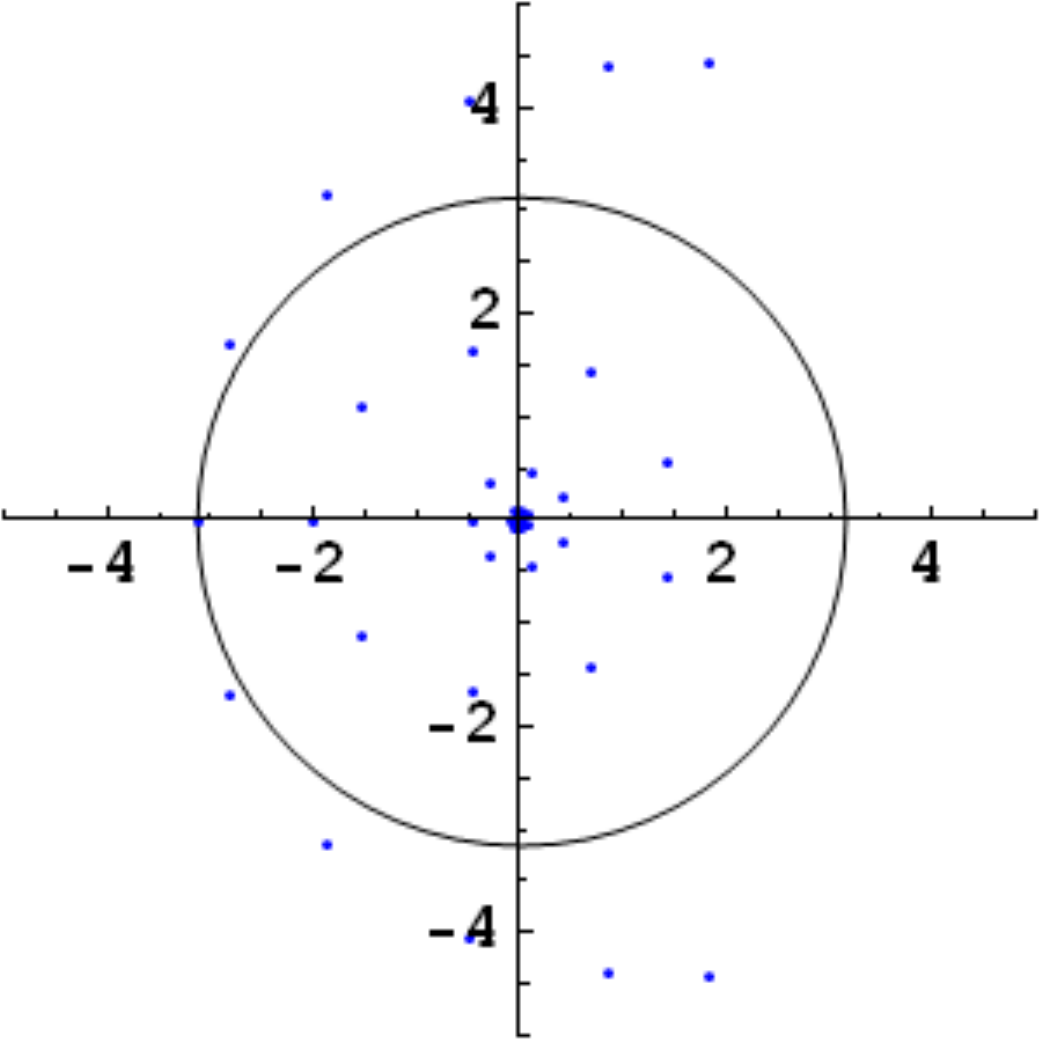}\\
\parbox{.47\textwidth}{\caption{$\lS_{7,24}(\zetastar)$}\label{pic724}}
\hfill\parbox{.47\textwidth}{\caption{$\lS_{7,48}(\zetastar)$}\label{pic748}}
\end{figure}
\begin{figure}[p]
\includegraphics[width=.47\textwidth]{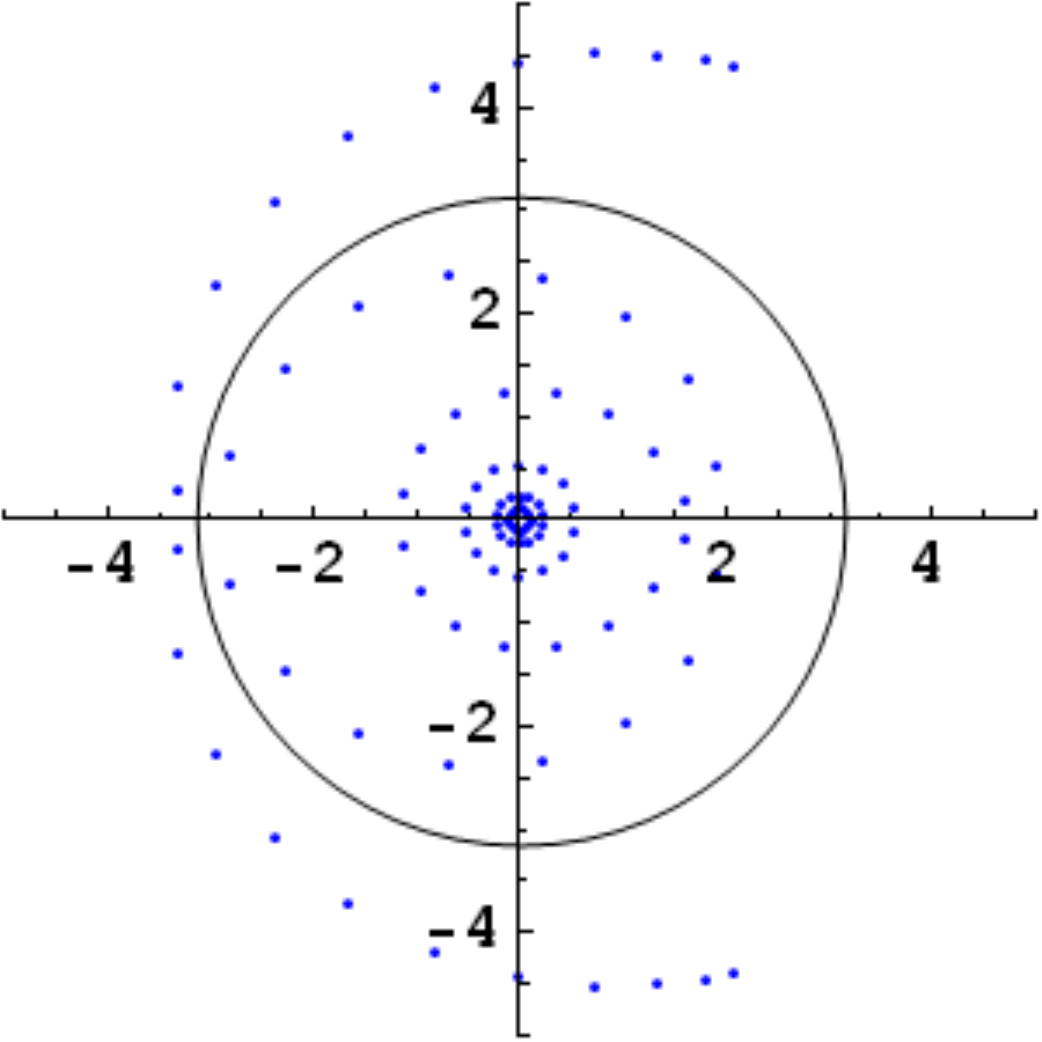}\hfill
\includegraphics[width=.47\textwidth]{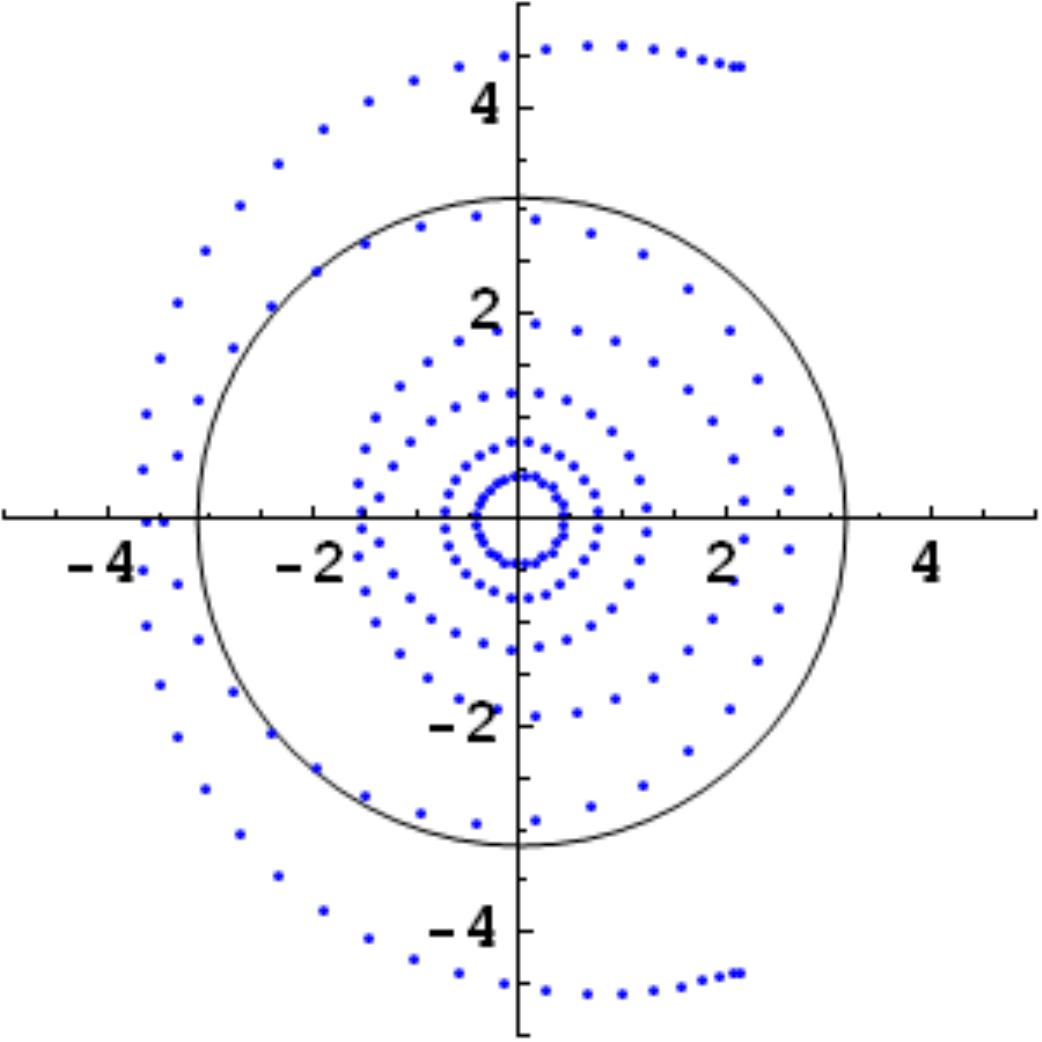}\\
\parbox{.47\textwidth}{\caption{$\lS_{7,96}(\zetastar)$}\label{pic796}}\hfill
\parbox{.47\textwidth}{\caption{$\lS_{7,192}(\zetastar)$}\label{pic7192}}
\end{figure}

The fifth orbit in the target appears in the spectra $\lS_{7,m}(\zetastar)$.
 Figures \ref{pic724}--\ref{pic7192} show    spectra $\lS_{7,m}(\zetastar)$
for $m=24,\ 48,\ 96, \ 192$
  respectively.
%Figure \ref{7192} shows  the union $\cup_{m=1}^{192}\lS_{7,m}(\zetastar)$.
An
animation showing the  $\lS_{7,1}(\zetastar)$, $\lS_{7,2}(\zetastar)$, \dots in succession can
be downloaded from \cite{hidden}.

The above pictures don't show the arrows  in any of
$\lS_{5,m}(\zetastar)$, $\lS_{6,m}(\zetastar)$,
$\lS_{7,m}(\zetastar)$ for $m=24,\  48,\ 96,\ 192$  but probably
the arrows will appear
for sufficiently large $m$.

\subsection{More Conjectures}

The above pictures suggest the following conjectures which, in particular,
 generalize conjectures
$\mathrm{1A_1}$--$\mathrm{1F_1}$.

\

{\bf Conjecture 1A.} \emph{There are never multiple eigenvalues in $\lS_{l,m}(\zetastar)$.}

\

{\bf Conjecture 1B.} \emph{For all $l$, $\sup_m( \max (\Arr_{l,m}(\zetastar)))$
is  a positive number.}

\

{\bf Conjecture 1C.} \emph{For all $l$, $\inf_m( \min (\Arr_{l,m}(\zetastar)))$ is
a positive number.}

\

Accepting the same agreement about the largest real eigenvalue in $\lS_{l,m}(\zetastar)$
as was done above in the case of
$\lS_{1,m}(\zetastar)$, we can state the following two conjectures.

\

{\bf Conjecture 1D.} \emph{For all $l$ and $k$, the numbers
$\arr_{l,m}(\zetastar)=||\Arr_{l,m}(\zetastar)||$,
$\orb_{l,m,k}(\zetastar)=||\Orb_{l,m,k}(\zetastar)||$
$\bow_{l,m}(\zetastar)=||\Bow_{l,m}(\zetastar)||$
of eigenvalues belonging to
the arrow $\Arr_{l,m}(\zetastar)$,
the orbit $\Orb_{l,m,k}(\zetastar)$, and
the bow $\Bow_{l,m}(\zetastar)$ respectively
 don't decrease when $m$ increases.}

\

{\bf Conjecture 1E.} \emph{For all $l$, if
\begin{eqnarray}
\Arr_{l,m}(\zetastar)&=&
\{\lambda_{l,m,1}(\zetastar),\dots,\lambda_{l,m,\arr_{l,m}(\zetastar)}(\zetastar)\},\\
\Arr_{l,m+1}(\zetastar)&=&\{\lambda_{l,m+1,1}(\zetastar),\dots,
\lambda_{l,m+1,\arr_{l,m+1}(\zetastar)}(\zetastar)\}
\end{eqnarray}
and
\begin{eqnarray}
&\lambda_{l,m,1}(\zetastar)<\lambda_{l,m,2}(\zetastar)<\dots<
\lambda_{l,m,\arr_{l,m}(\zetastar)}(\zetastar),\\
&\lambda_{l,m+1,1}(\zetastar)<\lambda_{l,m+1,2}(\zetastar)<\dots<
\lambda_{l,m+1,\arr_{l,m+1}(\zetastar)}(\zetastar)
\end{eqnarray}
then
\begin{equation}
\lambda_{l,m+1,1}(\zetastar)<\lambda_{l,m,1}(\zetastar),
\dots,
\lambda_{l,m+1,\arr_{l,m}(\zetastar)}(\zetastar)<\lambda_{l,m,\arr_{l,m}(\zetastar)}(\zetastar).
\end{equation}
}

\

{\bf Conjecture 1F.}
\emph{For given  $l$ and $m$,  assign the weight $\frac{1}{m}$ to each point
$\lambda_{l,m,1}(\zetastar),\lambda_{l,m,2}(\zetastar),\dots,\lambda_{l,m,m}(\zetastar)$
and denote by $\lambda_{l,m}(\zetastar)$ the corresponding discrete
measure.  Then
\begin{itemize}
\item[$\mathbf{1F'.}$] for $m\rightarrow \infty$
there exists a limiting continuous measure $\lambda_{l}^{\zetastar}(w)$
concentrated on a ``limiting bow'',  a ``limiting arrow'',
and a ``limiting target'' consisting of a number of
``limiting orbits'';
 \item[$\mathbf{1F''.}$]
$\int \log(w)\myd\lambda_{l}^{\zetastar}(w)=\log(\W_l)$.
\end{itemize}
}

\

Clearly, Conjecture $1\mathrm{F}$ implies the Riemann Hypothesis.

An precise definition of the bow $\Bow_{l,m}$ and
of individual orbits $\Orb_{l,m,k}$ constituting the
target $\Targ_{l,m}$ is a subtle matter because of the
rendezvous. The following conjectures implicitly
presuppose that \emph{It is possible to split $\lS_{l,m}(\zetastar)$
into $\Arr_{l,m}(\zetastar)$, $\Orb_{l,m,k}(\zetastar)$,
and $\Bow_{l,m}(\zetastar)$ in such a manner that...}.

\

{\bf Conjecture $\mathbf{1G}$.} \emph{For every $l$ and $k$  the number of
eigenvalues
in $\Arr_{l,m}(\zetastar)$, $\Orb_{l,m,k}(\zetastar)$,
and $\Bow_{l,m}(\zetastar)$ doesn't decrease with the growth of $m$.}

\

{\bf Conjecture $\mathbf{1H}$.} \emph{For every $l>1$, $m$, and $k$ the
eigenvalues from the
orbit $\Orb_{l,m,k}(\zetastar)$ almost lie on a circle and  are
almost equidistributed on it.}

\

{\bf Conjecture $\mathbf{1I}$.} \emph{For every $l>1$ the limiting
target $\Targ_l(\zetastar)$ consists of some number $\targ_{l}(\zetastar)$ of
limiting orbits
$\Orb_{l,k}(\zetastar)$ and
\begin{itemize}
\item[$\mathbf{1I}'$] all limiting orbits are circles;
\item[$\mathbf{1I}''$]
on each limiting orbit  the limiting measure~$\lambda_l^{\zetastar}(w)$ is constant;
\item[$\mathbf{1I}'''$] the limiting
 orbits
$\Orb_{l,k}(\zetastar)$
can be numbered in such a way that for $k<\targ_{l}(\zetastar)$
the limiting orbit $\Orb_{l,k}(\zetastar)$
lies inside the limiting orbit $\Orb_{l,k+1}(\zetastar)$ and touches it
at one real point $\Rend_{l,k}(\zetastar)$ called the \emph{rendezvous-point};
the innermost limiting orbit has rendezvous point $\Rend_{l,0}(\zetastar)$
with the limiting arrow $\Arr_{l}(\zetastar)$, and the
outmost limiting orbit has rendezvous point $\Rend_{l,\targ_{l}(\zetastar)}(\zetastar)$
with the limiting bow $\Bow_{l}(\zetastar)$;
moreover, $\Rend_{l,k+1}(\zetastar)<\Rend_{l,k}(\zetastar)$
for even $k$  and
$\Rend_{l,k+1}(\zetastar)>\Rend_{l,k}(\zetastar)$ for odd $k$.
\end{itemize}}

\section*{Acknowledgement}

\

The author is very grateful to Martin Davis
for some help with the English.

\end{document}